
\documentstyle{amsppt}

\magnification=1200
\pagewidth{6 truein}
\pageheight{9 truein}

\def \a {\alpha}
\def \ad {{\text {ad}}}

\def \affg {\operatorname{Aff}(\frak g)}
\def \affgs {\operatorname{Aff}(\frak g, \sigma)}
\def \Aff{\operatorname{Aff}}
\def \andd {\quad\text{and}\quad}
\def \Aut{\operatorname{Aut}}
\def \bari {{\bar \imath}}
\def \barj {{\bar \jmath}}
\def \C {\Bbb C}

\def \dg{\dot \g}
\def \dh{\dot \h}

\def \eala {(\g,\h, (\cdot,\cdot))}

\def \fm {(\cdot,\cdot)}
\def \g {\frak g}
\def \gia {\g_{\bar \imath, \pi (\alpha)}}
\def \gjb {\g_{\bar \jmath, \pi (\beta)}}
\def \gke {\g_{\bar k, \pi ( \epsilon)}}
\def \gt {\tilde \frak g}
\def \h {\frak h}

\def \id{\text{id}}

\def \ital#1{{\it #1\/}}  
\def \k {\frak k}
\def \la {\langle}
\def \lsa {\ell_{\sigma}(\alpha)}
\def \m {\frak m}
\def \minv{\frac 1m}  
\def \boldmu {{\pmb \mu}}

\def \ot {\otimes}
\def \p {{\pmb p}}
\def \R {\Bbb R}
\def \ra {\rangle}
\def \Rt {\widetilde R}
\def \Rx {R^{\times}}
\def \s {\sigma}
\def \ss {\frak s}
\def \suchthat{\,\mid \, }  
\def \t {\tau}
\def \ta {\tilde \alpha}
\def \tc {\tilde \frak c}
\def \tg {\tilde \frak g}
\def \th {\tilde \frak h}
\def \tR {\tilde R}
\def \ts {\tilde \frak s}
\def \tV {\tilde V}
\def \Z {\Bbb Z}

\def \zn {\Bbb Z^{\nu}}

\def \loopalg {\operatorname{L}}
\def \loopg {\operatorname{L}(\frak g)}
\def \loopgs {\operatorname{L}(\frak g, \sigma)}
\def \ident {\operatorname{id}}
\def \A {\Cal A}
\def \K {\Cal K}
\def \q {{\bold q}}
\def \sll{\frak{sl}}
\def \tr{\operatorname{tr}}
\def \e{\varepsilon}
\def \inv{\bar{}}
\def \Om{\Omega}
\def \Dl{\Delta}
\def \Dls{\Delta_{\text{v}}}
\def \Dlns{\Delta_{\text{nv}}}

\def\noi{\noindent}

\topmatter
\title
 Covering Algebras I:  
 Extended Affine Lie Algebras 
\endtitle
\author Bruce Allison\\
Stephen Berman\\
Arturo Pianzola\\
\endauthor
\address
Department of Mathematics, University of
Alberta,
Edmonton, Alberta, Canada T6G 2G1
\endaddress
\address
Department of Mathematics and Statistics,
University of
Saskatchewan, Saskatoon, Saskatchewan, Canada S7N 5E6
\endaddress
\address
Department of Mathematics, University of
Alberta,
Edmonton, Alberta, Canada T6G 2G1
\endaddress
\thanks The authors gratefully acknowledge the support
of the
Natural Sciences and Engineering Research Council of Canada.
\endthanks
\rightheadtext{Covering algebras I: EALA's}
\leftheadtext{}
\endtopmatter

\document

This is the first, of what will be a sequence of three papers, dealing 
with a  generalization of certain parts of the beautiful 
work of V. Kac on finite order automorphisms of  finite 
dimensional complex simple Lie algebras. Recall 
that Kac (see [K2, Chapter 8] and [H, \S X.5]) built a Lie algebra from a pair 
$(\g, \s)$ comprised  of a finite order automorphism $\s$ of a
 finite dimensional   simple Lie 
algebra $\g$ over $\C$ algebra as follows. First from $\s$  he 
obtains the 
eigenspaces
$$\g_{\bari}=\{x \in \g \suchthat \s(x) = \zeta^i x \}, $$
where $m$ is the order of $\s$, $\zeta=e^{2 \pi i/m},$ 
$i \in \Z$ and $i \to \bar{\imath}$ is the natural map of $\Z \to 
\Z_m$ (here $\Z_m$ denotes the integers modulo $m$). He then 
constructs the Lie algebra
$$ \affgs :=\big(\bigoplus_{i \in \Z}\g_{\bari} \otimes   t^i \big)\oplus 
\C c \oplus \C d $$   where  $c$ is 
central, $d=t\frac{d}{dt}$ is the degree derivation so that $[d, 
x\otimes t^p]=px \otimes t^p, $ and 
$$[x \otimes t^p, y \otimes t^q]=[x,y] \otimes t^{p+q}+ 
p(x,y)\delta_{p+q,0}c$$
Here $x \in \g_{\bar p}$, $y \in \g_{\bar q} \text{ and } p,q \in \Z,$ 
while $(\,,\,)$ denotes  the Killing form of $\g.$ It is known that all affine Kac-Moody Lie
 algebras arise 
this way. When $\s$ is the identity,  one gets the 
untwisted affine algebra corresponding to $\g$ (see [K1] and [M]), 
whereas the other graph automorphisms, if they exist, 
lead to the twisted versions of these 
algebras (see [K1]).   By taking into account the conjugacy theorem of 
Peterson and Kac [PK],  one can 
say that in fact the above construction establishes a bijective correspondence 
between affine 
Kac-Moody Lie algebras and
conjugacy classes of graph automorphisms of the finite dimensional 
simple Lie algebras. Furthermore, the  complex 
Lie algebra $ \affgs$
depends only on the outer part of $\s.$ Thus the $ \affgs$'s, 
which at first 
appear to make up a very large class of Lie algebras, give up to isomorphism 
exactly the affine Lie algebras.  Kac goes on to 
use this information to 
get a classification of all finite order 
automorphisms of $\g.$  

Our goal in this sequence of 
three papers is to get, as complete as is possible, a picture of what 
all of the algebras $\affgs$ look like for other classes of Lie 
algebras $\g$ as explained below. But first a word on terminology. At 
different stages of our work we will find it convenient, and necessary, to also 
work with several other 
versions of the algebras $\affgs$. Thus we will work with the ``loop version'' 
with no central element and no derivation added. It is given by
$$\loopgs= \bigoplus_{i \in \Z} g_{\bari} \otimes t^i$$
with multiplication coordinate wise. We will also need to study the 
derived algebra of $\affgs.$ Of course we need to carefully 
distinguish these and will do so when necessary. However, we also 
think of them all simply as {\it covering algebras} of our original 
algebra $\g$ relative to the automorphism $\s.$ Thus, we use the 
term covering algebra loosely and will develop other specific names 
for our various algebras when necessary. For example,  we will call 
the algebra $\affgs$ {\it the affinization of} $\g$ {\it relative to the 
automorphism} $\s.$

What follows is a very brief account of the topics discussed in
the three papers.
 
\medskip
\noi {1.} {\it Extended affine Lie algebras (EALA's)}. These algebras are natural 
generalizations of affine Kac-Moody Lie algebras, and are the subject 
of study in this our first paper. EALA's come equipped with analogues 
of Cartan subalgebras, root systems, invariant forms, etc. The role played by the null roots 
in the affine Kac-Moody case, is  that of the so called isotropic 
roots. These generate a lattice whose rank, is referred to as the 
nullity of the EALA. In fact, the result of Kac's theorem  in the 
language of EALA's reads as follows:  If
$\g$  is a tame EALA of nullity zero then $ \affgs$ is a tame EALA 
of nullity one and moreover, all such algebras arise in this way. 
When phrased  this way it becomes quite natural to ask what happens 
in the case of EALA's of higher nullity. In the first paper
we look at when the 
affinization $\affgs$ of an EALA is itself an EALA, and determine the 
relationship between  the corresponding root systems.
\medskip
\noi {2.} {\it Symmetrizable Kac-Moody Lie algebras}. The second 
paper begins with a general construction that provides a cohomological description of covering 
algebras by means of forms (a non-abelian $H^{1}$). The main thrust 
of the work lies in trying to decide if in the symmetrizable Kac-Moody 
case, the algebras $\affgs$ depend only on the outer part of $\s.$ 
(As we have seen, this is exactly what happens if $\g$ is finite 
dimensional.) To study this question, we use the Gantmacher-like description of 
automorphisms provided by [KW].
\medskip
\noi {3.} {\it Affine Kac-Moody Lie algebras.} In the final paper we 
study in detail the finer structure of the EALA's that are obtained 
as affinizations of the affine Kac-Moody Lie algebras.
\nopagebreak
\medskip
\nopagebreak
Each of the papers will have its own introduction.

\head {\S 1 Introduction: Basics on EALA's and an outline of results}\endhead

In this introduction to the present paper, we begin by recalling the definition of an extended affine Lie algebra (EALA for short) and some of the basic properties of these
algebras.  We will conclude with a short outline of the main results of the paper.

EALA's were introduced by R.~H\o egh-Krohn and B.~Torresani in [H-KT].
Many of the basic facts about EALA's and their root systems were
proved in [AABGP].  The reader can consult that reference
for any results stated in this section without proof.

Throughout this paper we will work with Lie algebras over the field of complex
numbers $\Bbb C$. The basic
definition of an EALA is broken down into a sequence of axioms EA1--EA4, 
EA5a, and EA5b.
 
Let $\g$ be a Lie algebra over $\C$. We assume first of all that $\g$ satisfies 
the following axioms EA1  and EA2.

\smallskip
{\bf EA1.} $\g$ has a non-degenerate invariant symmetric bilinear form
denoted by
$$(\cdot ,\cdot):\g \times \g \to \C.$$

{\bf EA2.} $\g$ has a nonzero finite dimensional abelian subalgebra
$\h$ such that $\ad_{\g}h$ is diagonalizable for all $h \in \h$ and
such that $\h$ equals its own centralizer, $C_{\g}(\h)$, in $\g$.

   One lets $\h^*$ denote the dual space of $\h$ and for  $\a \in
\h^*$ we let
$$\g_{\a}=\{ x \in \g \suchthat [h,x]=\a(h)x \text { for all
} h \in \h \}.$$
Then we have that
$\g = \bigoplus_{\a \in \h^*}\g_{\a}$ and so since 
$C_{\g}(\h)=
\g_0$ by EA2,
we obtain that $\h = \g_0.$

   We next define the \ital{root system} $R$ of $\g$ relative to $\h$ 
by
saying
$$R=\{\a \in \h^* \suchthat\g_\a \neq 0 \}.$$
Notice that $0 \in R$ and that $R$ is an extended affine root system
(EARS for short) in the sense of [AABGP]. As usual one finds that
$(\g_\a,\g_{\beta})=0$ unless $\a + \beta =0$. Thus, $-R=R$ and we
also have that $(\g_\a,\g_{\beta})=0$ unless $\a +\beta =0$. In
particular, the form is nondegenerate when restricted to $\h \times
\h$ and so this allows us to transfer the form to $\h^*$ as follows.
For each $\a \in \h^*$ we let $t_\a$ be the unique element in $\h$
satisfying
$(t_\a,h)=\a (h)$ for all $h \in \h.$ Then for $\a, \beta \in \h^*$ we
let $(\a, \beta)$ be defined by the equation
$$(\a, \beta)=(t_\a,t_{\beta}). \tag 1.1 $$

   We now have a nondegenerate form on $\h^*$ and so can speak of
isotropic and nonisotropic roots. We let
$R^0$ be the set of isotropic roots and let $R^{\times}$ be the
nonisotropic roots so that we have the disjoint union
$$R=R^0 \cup R^{\times}.$$

   We can now state the remaining axioms.

\smallskip
{\bf EA3.} For any $\a \in R^{\times}$ and any  $x \in \g_\a$ the
transformation $ \ad_\g x$  is a locally nilpotent
on $\g$.

{\bf EA4.} $R$ is a discrete subspace of $\h^*$.

{\bf EA5a.} $R^{\times}$ cannot be decomposed into a union
$R^{\times}=R_1 \cup R_2$ where $R_1$ and
$R_2$ are nonempty orthogonal subsets of $R^{\times}$.

{\bf EA5b.} For any $\delta \in R^0$ there is some $\a \in R^{\times}$
such that $\a + \delta \in R.$

\definition{Definition 1.2} A triple $(\g,\h,\fm)$
consisting of a Lie algebra $\g$, a subalgebra $\h$ and a bilinear
form $\fm$  satisfying EA1-EA4, EA5a, and EA5b is called an
\ital{extended affine Lie algebra} or EALA for short. \enddefinition

  We often will abuse notation and simply say  ``Let $\g$ be an EALA''
but the reader should always recognize that we have in mind a fixed
triple $\eala$.

Until further notice we 
let $\g$ denote an EALA.  We next recall
some properties of such algebras. 

There are always nonisotropic roots
as is shown by applying EA5b to the root $0$ (which is a root by EA2).
We have $(R^0,R^{\times})=\{ 0 \}.$ Let $V$ be the real span of $R$.
Then one knows (see Theorem 2.14
of [AABGP]) that the form on $\g$ can be scaled in such a way that
$(\a, \beta) \in \R $ for all $\a, \beta \in R$ and the
form is positive semidefinite on $V$.
Let $V^0$ be the radical of $V$.  Then
the integer $\nu=\dim_\R V^0$ is called
the \ital{nullity} of the EALA $\g$. We have that $V^0$ is the real 
span of $R^0$.
We denote the natural map from $V$ to $\bar V = V/V^0$
by $x \mapsto \bar x$, and we let $\bar R$ be the image of $R$ in $\bar V$. Then $\bar R$ is a
finite irreducible  root system in the space $\bar V$ where we use the
positive definite form on $\bar V$ induced from the semidefinite form on $V$.
($\bar R$ contains 0 and is possibly nonreduced.)
The \ital{type} of $\g$ is by definition the type of the 
root system $\bar R.$

It is also often useful to know that $\dim \g_\a =1 \text{ for all } \a \in \Rx$.

\definition{Definition 1.3} The \ital{core} of the EALA $\g$, denoted
by $\g_c$, is the subalgebra of $\g$ generated by the root spaces
$\g_\a$ for nonisotropic roots $\a \in \Rx.$ \enddefinition

It is easy to see that the core $\g_c$ of $\g$ is an ideal of $\g$.

\definition {Definition 1.4} We say the EALA $\g$ is \ital{tame} if 
the centralizer 
$C_{\g}(\g_c)$ of $\g_c$ in $\g$
is contained in $\g_c$.  Equivalently $\g$ is tame if  
$C_{\g}(\g_c)$ equals the
center $Z(\g_c)$ of $\g_c$.
\enddefinition

We are now in a position to briefly describe the contents of this paper.

In Section 2  we give the general definition of $\loopgs$ and $\affgs$.
For the definition of $\loopgs$ all that is required is a Lie algebra
$\g$ and a finite order automorphism $\s$ of $\g$.
To define $\affgs$ one requires in addition a 
nondegenerate invariant symmetric bilinear form on $\g$. 
As already mentioned, 
we call $\affgs$ the affinization of $\g$ relative to $\s$. 
$\affgs$ can  be realized as the subalgebra of 
fixed points of an automorphism of the Lie algebra 
$\operatorname{Aff}(\frak g, \ident_\g)$ (which is just the usual 
affinization of $\g.$) \ It is always the case that $\affgs$ has a nondegenerate 
symmetric bilinear form which is invariant.

Section 3 contains the main results of this paper.
We work with an
EALA $\eala$ and an automorphism $\s$ of $\g$ which satisfies the 
following four properties: 

{\bf A1.} $\s^m=1.$

{\bf A2.} $\s(\h)=\h.$

{\bf A3.} $(\s(x),\s(y))=(x,y)$ for all $x,y \in \g.$

{\bf A4.} The centralizer of $\h^{\s}$ in $\g^{\s}$ equals
$\h^{\s}.$

\noindent
Here we have let $\h^{\s}$ (respectively $\g^{\s}$) denote 
the fixed points of $\s$ in $\h$ (respectively~$\g$).

Assuming A1--A4, $\affgs$ has a natural choice
of a finite dimensional abelian ad-diagonalizable subalgebra,
namely $(\h^\s\otimes 1)\oplus \C c \oplus \C d$.
Also
$\affgs$ has a natural nondegenerate invariant symmetric bilinear form
$\fm$ which is the  form defined in section 2.
We
investigate whether or not the triple $(\affgs,
(\h^\s\otimes 1)\oplus \C c \oplus \C d,\fm)$
is an EALA by investigating whether axioms EA1, EA2, EA3, EA4, EA5a and EA5b hold. 
We find that the first five of these always hold.
In fact we see, in Proposition 3.25, that A4 is precisely
the assumption needed in order to obtain EA2.

Dealing with axiom EA5b, as well as 
whether or not $\affgs$ is tame,  is less straightforward. Our main result, 
Theorem 3.63, states that if $\g$ is tame  then 
either $\affgs$ has no nonisotropic roots or is again a tame EALA. 
It also describes how to tell which of these (mutually exclusive)
alternatives holds
using only information about the root system $R$ of $\g$ and the 
action of $\s$ on $R$. Thus, we give a very general procedure 
for producing new tame EALA's from old ones. 
This is accomplished without using the classification of tame
EALA's that is currently being developed (see [BGK], [BGKN], [Y] and [AG]). 

In the case when $\affgs$ is a tame EALA, 
we also obtain in Section 3 a description of the root system, 
core, type and nullity of $\affgs$
in terms of the corresponding objects for $\g$.  We further 
show that if $\g$ is non-degenerate (the definition of this term 
for EALA's is recalled when we need it), then  $\affgs$ is also nondegenerate.

As a byproduct of our investigation in Section 3,
we notice a general fact about the axioms of an EALA.
Namely EA5b 
follows from the other axioms of an EALA together with tameness. This 
may be of some independent interest.

Finally, in Section 4 we 
present 3 examples of how to use our main result to obtain EALA's
using affinization. 
In particular, in the third example we look at the 
important case when $\g$ is taken to be an affine Kac-Moody Lie 
algebra and $\s$ is a diagram automorphism. 
A1, A2 and A3 hold by definition of $\s$ and we show
A4 also holds in all 
cases.  
We find further that $\affgs$ has a nonisotropic root 
if and only if the diagram automorphism is not 
transitive on the nodes of the diagram.  
It follows from our results that
$\affgs$ is a tame EALA of nullity 2 in all cases but one.
The one exception is when $\g$ is of type 
$A_l^{(1)}$ and $\s$ is a  graph automorphism which is transitive 
on the set of nodes of the diagram. 
In the present paper we 
leave our investigation of covering algebras of affine algebras
here even though we fully  realize there is 
the natural question of identifying these algebras. This question will be 
investigated in the third paper in this series.
We note also that covering algebras of affine algebras 
have already been studied in 
[H-KT], [W] and [Po].

\head 
\S 2 The Lie algebras  $\loopgs$ and  $\affgs$
\endhead

In this section, we give the general definitions
of the Lie algebras $\loopgs$ and $\affgs$. 
 
Throughout the section we assume that
$\g$ is a Lie algebra over $\C$. (Although we work here with a complex Lie algebra,
the reader will notice that all that is really required in this section is that the
base field contain a primitive $m^{\text{th}}$ root of unity, where $m$ is a period for the automorphisms being considered.)

\definition{Definition 2.1} The {\it loop algebra\/} of $\g$
is the Lie algebra
$\loopg = \g \otimes \C[t,t^{-1}].$
\enddefinition

Suppose next that $m$ is a positive integer and $\s$ is an automorphism
of $\g$ of period m.  
Notice
that we do not require that $m$ is the actual
order of $\s$ but only a period of $\s$.
To define $\loopgs$, we need some notation.
Let $$\Z_m= \Z/m\Z$$
be the group of integers mod $m$
and let $i \to \bari$ denote
the natural homomorphism of $\Z$ onto
$\Z_m$.   Finally, let $\zeta= e^{2 \pi i/m}$.
Then $\g$ can be 
decomposed into eigenspaces for $\s$ as 
$$\g=  \bigoplus_{i=0}^{m-1} \g_{\bari},\quad 
\text{where }
g_{\bari}= \{ x \in \g \suchthat\s (x)= \zeta ^i x \}.$$ 
Observe that $\g_{\bar 0}$ is the subalgebra 
$\g^{\s}$ of fixed points of $\s$ in $\g$.

Now $\s$ extends to an automorphism  (also
denoted by $\s$)  of $\loopg$ defined by
$$\s (x \otimes t^i)=\zeta^{-i}\s (x)\otimes t^i.$$
The subalgebra $\loopg^{\s}$ of fixed points of
$\s$ in $\loopg$ is equal to
$\bigoplus_{i \in \Z}\g_{\bari} \otimes   t^i$.

\definition {Definition 2.2} Assume that $\s$ is an
automorphism of $\g$ such that $\s^m =1$. Let
$\loopalg(\g,\s,m)$ be the Lie algebra
$\bigoplus_{i \in \Z}\g_{\bari } \otimes   t^i$. 
It is easy to see that $\loopalg(\g,\s,m) \cong \loopalg(\g,\s ,n)$
if we have $\s^m =1 =\s^n$ and so, up to isomorphism, this algebra does not depend on the period $m$. Thus, we drop the dependence on $m$ and just write $\loopgs$ for
$\loopalg(\g,\s,m)$. So we have 
$$\loopgs =\loopg^{\s}= \bigoplus_{i \in \Z}\g_{\bari } \otimes   t^i$$
and $\loopgs$ is a subalgebra of the loop algebra $\loopg =
\loopalg(\g,\ident_{\g})$.   We
call $\loopgs$ the \ital{loop algebra of $\g$ relative to
$\s$}.
\enddefinition

If $0\ne a\in \C$, then the map $x\otimes t^k  \mapsto a^kx$ is a Lie algebra homomorphism
of $\loopg$ onto $\g$.  The restriction of this homomorphism maps $\loopgs$ onto $\g$.
This is the reason why $\loopgs$ (along with other related algebras) is regarded
as a covering algebra of~$\g$ (see [H, \S X.5]).

To discuss affinization, we assume that $\g$ is a Lie algebra over $\C$ with a
nondegenerate invariant symmetric bilinear form
$$(\cdot ,\cdot):\g \times \g \to \C.$$

\definition {Definition 2.3} 
The \ital{affinization} of $\g$ is defined as
$$\affg = \left(\g \otimes \C[t,t^{-1}]\right) \oplus \C c \oplus \C d.$$
Multiplication on $\affg$ is defined by
$$[x \otimes t^i +r_1c +r_2d, y \otimes t^j +s_1c +s_2 d]=
[x,y]\otimes t^{i+j}+j r_2y \otimes t^j-i s_2x \otimes t^i +
i\delta_{i+j,0} (x,y)c,$$
where $x,y \in \g,i,j \in \Z, r_1,r_2,s_1,s_2 \in \C.$ One easily checks that $\affg$ is
a Lie algebra.  We extend
the form on $\g$ to one on $\affg$ by decreeing that
$$ (x \otimes t^i +r_1c +r_2d,y \otimes t^j +s_1c +s_2 d)=
\delta_{i+j,0}(x,y)+r_1 s_2 +r_2 s_1, \tag 2.4 $$
and see immediately that this gives a nondegenerate invariant
symmetric bilinear form on $\affg$.
\enddefinition

Note that the above construction coincides with the usual affinization
in the case when $\g$ is a finite dimensional simple Lie algebra.
Also, $c$ is central in $\affg $ and $d$ is just the usual degree
derivation of $\C[t,t^{-1}]$ lifted to this affinization.  

To define relative affinization, we need the following lemma.

\proclaim {Lemma 2.5} Let  $\s$ be an automorphism
of $\g$ such that $\s^m =1$ and  $(\s x, \s y)=(x,y)$ for
all $x,y \in \g.$   Then $\s$ extends to an automorphism  (also
denoted by $\s$)  of $\affg$ defined by
$$\s (x \otimes t^i +rc +sd)=\zeta^{-i}\s (x) \otimes t^i
+rc+sd.\tag 2.6 $$
This extension preserves the extended form 2.4 on $\affg$, has
period $m$ and  fixes $c$ and $d$.   The subalgebra $\affg^{\s}$
of fixed points $\s$ in $\affg$ is given by
$$\affg^{\s} = \big(\bigoplus_{i \in \Z}\g_{\bari} \otimes   t^i\big)+\C
c+\C d.$$
Consequently, $(\bigoplus_{i \in \Z}\g_{\bari} 
\otimes   
t^i)
+\C c+\C d$ is a subalgebra of $\affg$ and the bilinear form on
$\affg$ is nondegenerate when restricted to this subalgebra.
\endproclaim

\demo{Proof}  All of the statements here are easily checked.
\qed
\enddemo

\definition {Definition 2.7} Assume that $\s$ is an
automorphism of $\g$ such that $\s^m =1$ and that
$(\s x, \s y)=(x,y)$ for all $x,y \in \g$. Let
$\Aff(\g,\s,m)$ be the Lie algebra
$(\bigoplus_{i \in \Z}\g_{\bari } \otimes   t^i)+ \C c+\C d$. 
If we have $\s^m =1 =\s^n$, it is easy to see that there
is a form preserving isomorphism of 
$\Aff(\g,\s,m)$ onto $\Aff(\g,\s,n)$.
So $\affgs$ does not depend on the period $m$,
and we just write $\affgs$ for
$\Aff(\g,\s,m)$.
Thus we have 
$$\affgs =\affg^{\s}=(\bigoplus_{i \in \Z}\g_{\bari } \otimes   t^i)+ \C c+\C d$$
and $\affgs$ is a subalgebra of the affinization $\affg =
\Aff(\g,\ident_{\g}).$  The form 2.4 restricted to $\affgs$ is a nondegenerate
invariant symmetric bilinear form on $\affgs$.
We call $\affgs$ the \ital{affinization of $\g$ relative to
$\s$}. 
\enddefinition

\head {\S 3 New EALA's from Old}\endhead

In this section we investigate the conditions required for
the affinization of an EALA relative to an automorphism of finite order
to  be an EALA.  We begin by specifying the assumptions 
that will be in force throughout the section.

\remark{\rm \bf Basic Assumption 3.1} Assume that $\eala$ is an 
EALA, $m
\in \Z$, $m \geq 1$ and $\s$ is automorphism of $\g$ which
satisfies
\smallskip
{\bf A1.} $\s^m=1.$

{\bf A2.} $\s(\h)=\h.$

{\bf A3.} $(\s(x),\s(y))=(x,y)$ for all $x,y \in \g.$
\smallskip

\noindent
Later we will add one more assumption, A4, on $\s$ (see Basic
Assumption 3.27).
\endremark

\medskip

Since we are mainly interested in the algebra $\affgs$, we make notation easier
by simply  writing
$$\tg = \affgs =\left(\g \otimes \C[t,t^{-1}]\right) \oplus \C c \oplus \C d.$$
We identify $\g=\g\otimes 1$ as a subalgebra of $\tg$, and we set
$$\th = \h^\s \oplus \C c \oplus \C d.$$
$\th$ is an abelian subalgebra of $\tg$.  Also,
we let $\fm$ denote the form on $\tg$ obtained by restricting the form 2.4
from $\affg$ to $\tg$.

Then specifically our problem in this section is to find conditions under
which the triple $(\tg,\th,\fm)$ is an EALA.

Using Lemma 2.5 we have the following:

\proclaim{Lemma 3.2} The form $\fm$ on $\tg$ is a nondegenerate invariant symmetric bilinear form.  Hence, $\tg$ satisfies EA1. \endproclaim

We next wish to consider EA2.  Thus, we need to 
look at the adjoint action of $\th$ on $\tg$.
To do this we will need to look at  the adjoint
action of $\h^\s$ on $\g$.
As in Section 1,  $R$ will denote the set of roots of $\g$ with
respect to $\h$. Also, for any vector space $X$ over $\C$, we let $X^*$ denote
the dual space of $X$ over $\C$.

 \definition {Definition 3.3} With notation as above we let $\pi :
\h^* \to (\h^{\s})^*$ be the linear map defined by saying
$\pi (\a)$ is just $\a$ restricted to $\h^{\s}$.
Symbolically, $\pi (\a)=\a|_{\h^{\s}}.$
\enddefinition

This mapping $\pi$ will play an important role in what follows. When
we consider $\h^{\s}$ acting on $\g$ via the adjoint action we see
that each element acts semisimply so that the set of weights for this
action is just $\pi (R).$ Furthermore, the eigenspace $\g_{\bar
\imath}$ for $\bari \in \Z_m$ is stabilized by $\h^{\s}$ so
we let
$\gia $ be the $\pi(\a)$ weight space of $\h^{\s}$ acting on
$\g_ {\bari}.$ That is
for all $\a \in R$ and all $\bari \in \Z_m$ we have
$$\gia =\{x \in \g_{\bari} \suchthat[h, x]= \a(h)x \text{
for all }h \in \h^{\s} \}. \tag 3.4 $$
We certainly have
$$ \g_{\bari}=\bigoplus _{\pi (\a) \in \pi (R)}\gia \text {
for any } i \in \Z.$$
Next, we let
$$R_{\bari} = \{ \a \in R\suchthat\gia \neq (0)\} \text {
for any } i \in \Z. \tag 3.5$$
Then we obtain
$$\g_{\bari}= \bigoplus_{\pi (\a) \in \pi(R_{\bari})}
\gia \text{ for } i \in \Z \text { and so } \tag 3.6 $$
$$\g= \bigoplus_{i=0}^{m-1} \bigoplus_{\pi (\a) \in
\pi(R_{\bari})} \gia. \tag3.7$$

Notice that we have if $\a, \beta \in R$  satisfy $\pi(\a)=\pi
(\beta)$ then $\a \in  R_{\bari }$ if and only if $\beta \in
R_{\bari}$ for any $i \in \Z.$ Moreover, for $\a \in R$ we
have that
$$\g_{\a} \subseteq \bigoplus_{i=0}^{m-1}\gia, \tag
3.8$$
which implies that $\a \in R_{\bari }$ for some $i \in \Z.$
Thus, we have the following union which is not necessarily disjoint:
$$R= \bigcup_{i=0}^{m-1}R_{\bari }.\tag3.9$$
We clearly have for $\a, \beta \in R, i,j \in \Z$ that either
$[\gia,\gjb]=(0)$ or there is some $\gamma \in R_{\bari
+\bar\jmath}$ with $\pi (\a) +\pi(\beta) =\pi (\gamma)$ and
$$[\gia,\gjb] \subseteq \g_{\bar\imath +\barj,\pi(\gamma)}. \tag
3.10$$

Next, we record how the spaces $\gia$ interact with our form. For $i,j
\in \Z$ and  $\a, \beta \in R$, if either
$\bari + \barj \neq 0$ or if $\pi(\a) + \pi(\beta)
\neq 0$ then we have
$$(\gia,\gjb)=(0).\tag3.11$$
Since the form is nondegenerate it follows that $-\a \in R_{-\bar
\imath}$ and the form pairs the spaces
$\gia$ and $\g_{-\bari,\pi(-\a)}$ in a nondegenerate
fashion.

Since we are assuming that $\s(\h) =\h$ we can decompose $\h$
relative to $\s.$ This gives us
$$ \h =\bigoplus_{i=0}^{m-1}\h_{\bari} \text{ and }
\h_{\bar 0}=\h^{\s}. \tag 3.12$$
Also, because $\s$ preserves the form we find that $(\h_{\bar
\imath},\h_{\barj}) =0$ whenever
 $\bari +\barj \neq 0.$ It follows that we have
$$(\h^{\s})^{\perp}=\bigoplus_{i=1}^{m-1}
\h_{\bari}. \tag 3.13$$  We therefore have
$$\h =\h^{\s} \oplus (\h^{\s})^{\perp} \quad\text{and
correspondingly}\quad
\h^* =(\h^{\s})^*\oplus ((\h^{\s})^{\perp})^*,$$
where we are identifying $(\h^{\s})^*$ and $((\h^{\s})
^{\perp})^*$ inside of $\h^*$
by letting any element from $(\h^{\s})^*$ act as zero on
$(\h^{\s})^{\perp}$
and letting any element from $((\h^{\s})^{\perp})^*$ act as zero
on
$\h^{\s}$.  With these identifications, note that $\pi :  \h^* \to
(\h^\s)^*$
is just the projection onto the first factor in the sum $\h^* =
(\h^{\s})^*\oplus ((\h^{\s})^{\perp})^*$.

 We let $\s$ acts on $\h^*$ by saying that
$$\s (\a)( h)=\a (\s^{-1}(h)) \text{ for all } h \in
\h. \tag 3.14 $$
This allows us to consider $(\h^*)^{\s}$, the fixed points of
$\s$ in $\h^*.$  It is not hard to see that $(\h^*)^{\s}=
(\h^{\s})^*$ and that $((\h^{\s})^{\perp})^*=((\h^*)^{\s})
^{\perp}$ and that this latter space is nothing but the sum of the
eigenspaces of $\s$ in $\h^*$ corresponding to eigenvalues
different from $1.$   For the convenience of the reader we next state
this as a lemma and sketch a proof.

\proclaim {Lemma 3.15}  With notation as above we have that 
$$(\h^*)^{\s}=(\h^{\s})^* \andd ((\h^{\s})^{\perp})^*=((\h^*)^{\s})^{\perp}$$ with this latter space being the sum of the  eigenspaces of $\s$ in $\h^*$
corresponding to eigenvalues different from $1.$
\endproclaim

\demo{Proof} We just prove the first statement as the second one is
similar. We let $\a \in (\h^{\s})^* \subseteq \h^*$ so that
$\a$ is zero on $(\h^{\s})^{\perp}.$ For $h \in \h^{\s}$
we have $\s (h)=h$ so that $\s^{-1}(h)=h$ and so by 3.14 we
get that $\s (\a)(h)=\a(h)$ so both
$\a$ and $\s(\a)$ agree on $\h^{\s}.$ If $h \in
(\h^{\s})^{\perp}$ then because $\s$ preserves the form so
does $\s^{-1}$ and so it follows that $\s^{-1}$ maps
$(\h^{\s})^{\perp}$ to itself. Thus $\a (\s^{-1}(h))=0$
and this implies that $\s(\a)(h)=0$ from which it follows that
both $\s(\a)$ and $\a$ agree on the space $(\h^{\s})
^{\perp}.$ This gives us $(\h^{\s})^* \subseteq (\h^*)^{\s}.$

   Conversely, let $\a \in (\h^*)^{\s}.$   Then $\a$ and
$\s(\a)$ agree on all of $\h .$ Next, recall 3.13 and choose
$\bari \in \Z_m \setminus \{ 0 \},$ as well as an element $h \in
\h_{\bari},$ so that $\s (h)=\zeta^i h$ and hence
$\s^{-1}(h)=\zeta^{-i}h.$ Thus, $\a(h)=\s(\a)(h)=
\zeta^{-i}\a(h)$ so since $\zeta^{-i} \neq 1$ we obtain that
$\a (h)=0. $Thus $\a$ and $\s(\a)$ are both zero on
$(\h^{\s})^{\perp}$ and so we now have $\a =\s(\a)$ is
in $(\h^{\s})^*$ as desired.
\qed
\enddemo

We are now ready to calculate the adjoint action of $\th$ on $\tg$.
Since
$$\th = \h^\s \oplus \C c \oplus \C d,$$
we can identify $(\h^{\s})
^*$  as a subspace of $\th^*$ so that the elements of $(\h^{\s})
^*$ act as zero on both $c$ and $d.$ We also define elements $\gamma$
and $\delta$ in
$\th^*$ by saying
$$\gamma(\h^{\s})=\{ 0\}=\delta(\h^{\s}),\quad \gamma(c)=1=
\delta(d),\quad \gamma(d)=0=\delta(c).\tag3.16$$
Thus, we have that
$$\th^*=(\h^{\s})^* \oplus \C \gamma \oplus \C \delta. \tag3.17$$

It is clear that $\th$ acts semisimply on $\tg$ via the adjoint
representation.  Thus,
if we set $\tg_{\tilde \a}=\{ x \in \tg \suchthat [h,x]=\tilde \a(h)x
\text { for all  } h \in \th \}$
for all $\tilde \a\in\th^*$, we
have
$\tg=\bigoplus_{\tilde \a \in \tg^*}\tg_{\tilde \a}.$
We set
$$\tR=\{\tilde \a \in \h^* \suchthat\g_{\tilde \a} \neq 0 \},$$
the set of roots of $\tg$ with respect to $\th$. Then,
$$\tg=\bigoplus_{\tilde \a \in \tR}\tg_{\tilde \a}. \tag 
3.18$$
Next 3.4 and 3.5 imply that if $R_{\bari}$ is nonempty then,
for any $i \in \Z$,  $\gia \otimes   t^i$ is decomposable into weight
spaces for $\th$ corresponding to
weights from $\pi (R_{\bari}) +i\delta.$  If $R_{\bari}$
is empty we take $\pi (R_{\bari}) +i\delta$ to be the empty set
for any $i \in \Z$.
We thus have
$$\tR= \bigcup_{i \in \Z}(\pi(R_{\bari})+i \delta) \subseteq
\th^*. \tag 3.19$$
Moreover,  if $i \in \Z$, $\a \in R_{\bari}$ and
$\pi(\a)+i \delta \neq 0$ then
$$\tg_{\pi(\a)+i \delta}=\gia \otimes   t^i. \tag3.20$$
On the other hand, we have that
$$\tg_0=(\g_{\bar 0,0 }\otimes   1 ) \oplus \C c \oplus \C d.
\tag3.21$$
But from 3.4 we see that $\g_{\bar 0,0 }$ is nothing but the
centralizer of $\h^{\s}$ in $\g^{\s}$. That is,
$$\tg_0=C_{\g^{\s}}(\h^{\s})\oplus \C c \oplus \C d.
\tag3.22$$

   Before stating our result about EA2, we need a result
about our map $\pi$ defined in 3.3. This is provided in the following
lemma.

\proclaim{Lemma 3.23} For $\a \in \h^*$ we have that $\pi(\a)=
\minv \sum_{i=0}^{m-1} \s^i(\a).$
\endproclaim

\demo{Proof}For $\a \in \h^*$ write $\a=\a_{\bar 0}+ \dots
+ \a_{\bar m -\bar 1}$ where $\a_{\bari} \in (\h^*)
_{\bari}.$ Then we have $\pi(\a)=\a_{\bar 0}.$ But
$\minv \sum_{i=0}^{m-1} \s^i(\a)=\a_{\bar 0}
+\minv\sum_{i=0}^{m-1}\s^{i}(\sum_{\barj \in
\Z_m \setminus \{ 0 \}}\a_{\barj}).$ Thus it is enough to
show that
$\sum_{i=0}^{m-1} \s^i(\a_{\barj})=0$ for
$\barj \in \Z_m \setminus \{ 0 \}.$ But we have
$\sum_{i=0}^{m-1} \s^i(\a_{\barj})=
\sum_{i=0}^{m-1}\zeta^{ij}\a_{\barj}
=  \frac{(\zeta^j)^m -1}{\zeta^j-1} \a_{\barj},$
which clearly
equals zero.
\qed
\enddemo

   We need one more definition before stating our result about EA2.
To prepare for this notice that
 for $\a \in R$ we have $\s(\g_{\a})=\g_{\s(\a)}$
and hence
$$\s(R)=R.$$

\definition{Definition 3.24} For $\a \in R$ we let $\lsa$ denote
the smallest positive integer satisfying
$\s^{\lsa}(\a) = \a.$ We call $\lsa$ the $\s$
\ital{length} of $\a$.
\enddefinition
For a root $\a \in R$ we have that $\a,\s( \a), \dots,
\s^{\lsa -1}(\a)$ are distinct. Moreover $\lsa$ divides the
period $m.$

\proclaim {Proposition 3.25} The following statements are equivalent.
\roster
\item"{(i)}" The centralizer of $\th$ in $\tg$ equals $\th.$ In other
words, EA2 holds for $\tg$ and $\th.$
\item"{(ii)}" The centralizer of $\h^{\s}$ in $\g^{\s}$ equals
$\h^{\s}.$
\item"{(iii)}" If $\a \in R \setminus \{ 0 \}$ then either \
$\pi(\a) \neq 0$ or \
$\{ x\in \g_a \suchthat \s^{\lsa}(x)  = x\} = \{0\}$.
\endroster
Moreover, if $m$ is prime these statements are equivalent to the
following.
\roster
\item"{(iv)}" If $\a \in R$ and if $\a \neq 0$ then
$\pi(\a) \neq 0.$
\endroster
\endproclaim

\demo{Proof}From 3.22 we see that (i) and (ii) are equivalent. We
first show that (ii) implies (iii). Thus, assume (ii) holds but that
(iii) fails. Then for some nonzero $\a \in R$ and
some nonzero $x \in \g_{\a}$ we have $\pi (\a)=0$ and
$\s^{\lsa}(x) =x .$ We write $\ell$ for $\lsa$ and let $y=
x+\s(x)+ \dots +\s^{\ell -1}(x).$ Now $y \neq 0$ as the roots
$\a, \s(\a),\dots ,\s^{\ell -1}(\a)$ are distinct
and we have $\s(y)=y$ so that
$y \in \g^{\s}.$ Since $\pi (\a)=0$ we have that
$\a(\h^{\s}) =\{0 \}$. Thus for $i=0,1, \dots , \ell-1$ we
have $(\s^i (\a))(\h^{\s})= \{ 0 \}$ and hence $[
\h^{\s},y]= \{ 0 \}.$ That is $y \in C_{\g^{\s}}(\h
^{\s}).$  Since $\a \neq 0$, we have $y \notin \h$ and so $y
\notin \h^{\s}$. This contradicts (ii).

   For the converse we assume now that (iii) holds but that (ii)
fails. Thus, there exists some $x \in C_{\g^{\s}}(\h^{\s})$ so
that $x \notin \h^{\s}.$  We can write $x$ as $x=
\sum_{\a \in R}x_{\a}$ where $x_{\a} \in \g_{\a}$
for all $\a \in R.$ Since $\s(x) =x$ it follows that
$\s^k( x_{\a})=x_{\s^k (\a)}$ for $\a \in R$, $k
\in \Z$, $k \geq 0.$ Thus $\s^{\lsa}(x_{\a})=x_{\a}$ for
$\a \in R.$ Also $[\h^{\s},x]=\{0\}$ so that
$[\h^{\s},x_{\a}]=\{0\}$ for $\a \in R.$ Thus
$\a(\h^{\s})x_{\a}= \{0\}$ for $\a \in R.$ Hence, if
$\a \in R$ and $x_{\a} \neq 0$ then $\a(\h^{\s})=\{0\}
$ and so $\pi (\a)=0.$ From (iii) it follows that $x_{\a}=0$
for all $\a \in R \setminus \{0\}$ and so we obtain that $x=x_0
\in \g_0 =\h.$  This gives us the contradiction that $x \in
\h^{\s}$ so establishes that (iii) implies (ii).

   Finally suppose that $m$ is prime. Clearly we only need to show
(iii) implies (iv) which we do by contradiction. Thus, assume there
exists some $\a \in R \setminus \{0\}$ such that $\pi(\a)=0.$
By Lemma 3.23 we have that $\s (\a) \neq \a$ and hence
$\lsa\neq 1.$ But $\lsa$ divides $m$ which is prime and so $\lsa =m$.
Hence
$\s^{\lsa}(x_{\a})=x_{\a}$ for all $x_{\a} \in
\g_{\a}$ and this contradicts (iii).
\qed
\enddemo

The conditions in the proposition automatically hold if for any
$\a \in R \setminus \{0\}$ we have $\pi (\a) \neq 0.$ Thus we
have the following corollary.

\proclaim{Corollary 3.26} If for any $\a \in R \setminus \{0\}$ we
have $\pi (\a) \neq 0 $ then $(\tg,\th,\tR)$ satisfies EA2.
\endproclaim

   Because of Proposition 3.25 we see that it is important for us that
the automorphism $\s$ of  $\g$ 
satisfies  $C_{\g^{\s}}(\h^{\s})=\h^{\s}.$  We will
require this as an assumption from now on.

\remark{\rm \bf Basic Assumption 3.27} We now assume that the 
automorphism
$\s$ of the EALA $\g$ satisfies A1, A2, A3 and
\smallskip
{\bf A4.} The centralizer of $\h^{\s}$ in $\g^{\s}$ equals
$\h^{\s}.$
\endremark

\medskip
\remark{\rm \bf Remark 3.28} Because of this basic assumption we see  
from
Lemma 3.2 and Proposition 3.25 that $(\tg,\th,\fm)$
satisfies
EA1 and EA2.
\endremark

\medskip
We now go on to consider the axioms EA3, EA4, EA5a and EA5b for $\tg$.

   Recall that $R^0$ and $R^ {\times}$ denote the set of isotropic and
nonisotropic roots respectively in R. Similarly we let $\tR^0$ and
$\tR^{\times}$ denote the set of isotropic and nonisotropic roots
respectively of $\tR.$ 

 From 3.16 and 3.17 we have that
$$(\gamma,\gamma)=0=(\delta,\delta) \andd (\gamma, \delta)=1. \tag 3.29
$$
This follows immediately from the way in which we transfer the form
from $\th$ to $\th^*.$ The following fact, which is easy to see, should be noted.
The restriction of the form on $\th$ to $\h^{\s}$ is
the same as the restriction of the form from $\h$ to $\h^{\s}.$ It
follows also that the restriction of the form from $\th^*$ to
$(\h^{\s})^*$ equals the restriction of the form from $\h^*$ 
to~$(\h^{\s})^*.$

   Now 3.19 says that $\tR= \bigcup_{i \in \Z}(\pi(R_{\bari})+i
\delta).$ Also for $i \in \Z$, $\a \in R_{\bari}$ we have $(\pi (\a )
+i \delta ,\pi (\a ) +i\delta )=
(\pi (\a ),\pi (\a))$ where now the right hand side can be
interpreted in $\h^*$ or in $\th^*.$ Thus we obtain that
$$\gather
\tR^{\times}= \bigcup_{i \in \Z} \{ \pi (\a )+ i \delta \suchthat
\a \in R_{\bari}, \ (\pi(\a), \pi (\a)) \neq 0 \}
\quad\text{and} \tag 3.30 \\
\tR^0=\bigcup_{i \in \Z} \{ \pi (\a ) +i \delta \suchthat\a
\in R_{\bari}, \ (\pi(\a ), \pi ( \a ))=0 \}. \tag
3.31\endgather$$

\proclaim  {Lemma 3.32} Let $\tilde \a \in \tR^{\times}$ and $x \in
\tg_{\tilde \a}.$ Then the operator
$\ad_{\tg}(x)$ is nilpotent. Thus $\tg$ satisfies EA3.
\endproclaim

\demo{Proof} We have that $\tilde \a =\pi(\a ) +i \delta$
where $i \in \Z$, $\a \in R_{\bari}$ and
$(\pi(\a ), \pi( \a )) \neq 0.$ Then $x \in \g_{\bari
,\pi (\a )} \otimes   t^i$ and we put
$h=t_{\pi(\a )} \in \h^{\s}$ (recall 1.1). Then 
$(\pi(\a ))(h)= ( \pi( \a ), \pi (\a ))$ and so $[h, x]=
( \pi( \a ), \pi (\a ))x.$ Hence $x$ is an eigenvector for
$\ad (h)$ corresponding to a non-zero eigenvalue. Thus, to show that
$\ad_{\tg}(x)$ is nilpotent it suffices to show that the
diagonalizable operator
$\ad_{\tg}(h)$ has only finitely many eigenvalues.

   By 3.19 the eigenvalues of $\ad_{\tg}(h)$ are the scalars of the
form $(\pi(\epsilon )+ j \delta)(h)$ where
$j \in \Z$, $\epsilon \in R_{\barj}.$ But since $h= t_{\pi
(\a )} \in \h^{\s}$ we have
$(\pi(\epsilon )+ j \delta)(h)=(\pi (\epsilon ))(h)=(\pi(\epsilon ),
\pi (\a ))=(\epsilon, \pi(\a )).$ By Lemma 3.23 we see this
equals $\minv(\epsilon, \sum_{i=0}^{m-1}\s^i(\a )).$
Thus, each eigenvalue of $\ad _{\tg}(h)$ is in the set $\minv(A +
\dots + A)$ ($m$ summands), where $A=\{ (\epsilon, \nu)
\suchthat\epsilon, \nu \in R\}.$ But we know that the set $A$ is
finite (see 2.11 and 2.14 in Chapter 1  and 2.10 in Chapter 2 of
[AABGP]).
\qed
\enddemo

   \proclaim {Lemma 3.33} $\tR$ is a closed and discrete subspace of
$\th^*.$ In particular, $\tg$ satisfies~EA4.
\endproclaim

\demo{Proof}We must show that any sequence in $\tR$ which converges to
an element of $\th^*$ is eventually constant. Looking at 3.19 we see
that such a sequence must eventually have the $i\delta$-part constant.
So it is enough for us to show that for fixed $i$ any sequence in
$\pi(R_{\bari})$ which converges to an element of $(\h^{\s})
^*$ is eventually constant. Thus, it suffices to show that any
sequence $\{ \pi(\a_n)\}_{1 \leq n < \infty}$ in $\pi(R)$ which
converges to an element of $\h^*$ is eventually constant. Now by Lemma
3.23 together with the fact that $\s(R)=R$ we get that
$\pi(R) \subseteq \minv(R+ \dots +R)\subseteq \minv\la R\ra$ where
$\la R \ra$ denotes the additive subgroup of $\h^*$ generated by $R.$
Thus $\pi(\a_n) \in \minv\la R \ra.$ But we know from 2.18 in Chapter II of [AABGP]
that $\la R \ra$ is
a lattice in the real span of $R.$ Thus so is $\minv\la R \ra.$ But
the real span of $R$ is closed in $\h^*$ and so our sequence $\{
\pi(\a)_n\}_{1 \leq n < \infty}$ converges to an element of this
span. Hence our sequence is eventually constant as desired.
\qed
\enddemo

We next want to investigate axiom EA5a.  In order to do this we will
need a lemma about finite irreducible root systems.  Since
we are using this term in a slightly nonstandard sense, we make our definition precise.
A \ital{finite irreducible root system} $\Om$ in a nonzero real Euclidean space $X$
is a finite spanning set for $X$ so that
$0\in \Om$, $2(\a,\beta)(\beta,\beta)^{-1}\in \Z$ 
for $\a$, $\beta$ in $\Om^\times$,
$\a - 2(\a,\beta)(\beta,\beta)^{-1}\beta\in \Om^\times$ 
for $\a$, $\beta$ in $\Om^\times$, and $\Om^\times$
cannot be written as the union of two nonempty orthogonal sets,
where $\Om^\times$ denotes the set of nonzero elements of $\Om$.

\proclaim{Lemma 3.34}  
Let $\Om$ be a finite irreducible root system in a nonzero real Euclidean space $X$.
Let $Y$ be a nonzero subspace of $X$, and let $p : X \to Y$ be the orthogonal
projection onto $Y$.  Then
\roster
\item"{(i)}" If $\a\in \Om^\times$,  there exists $\beta\in \Om^\times$ so that
$p(\beta)\ne 0$ and $(\a,\beta) \ne 0$.
\item"{(ii)}" $p(\Om^\times)\setminus\{0\}$ cannot be written as the union of
two nonempty orthogonal sets.
\endroster
\endproclaim

\demo{Proof} It will be convenient to introduce some temporary notation and terminology. 
Let $\Dl = \Om^\times$ (and so $\Dl$ is a root system in the usual sense).
If
$\a\in \Dl$, we say that $\a$ is \ital{visible} if $p(\a) \ne 0$.  Also, if $\a\in \Dl$,
we say that $\a$ is \ital{nearly visible} if there exists $\beta\in \Dl$ so that $\beta$
is visible and $(\a,\beta)\ne 0$.
We let $\Dls$ (resp.~$\Dlns$) be the set of visible (resp.~nearly visible) roots of $\Dl$.
Note that $\Dls$ is a subset of $\Dlns$.  Also, since $\Dl$ spans $X$ and since $Y\ne \{0\}$, we
have $p(\Dl) \ne \{0\}$ and hence $\Dls$ is nonempty.

We now claim that
$$\a\in \Dlns\setminus\Dls 
\implies 
\foldedtext\foldedwidth{2.5truein}
{There exist  $\beta, \gamma\in \Dls$
such that $\alpha = \gamma - \beta$
and $p(\gamma) = p(\beta).$}
\tag 3.35$$
Indeed, suppose that $\a\in \Dlns\setminus\Dls$. Then, $(\a,\beta)\ne 0$ for some
$\beta\in \Dls$.  In that case, $\a \ne \pm \beta$ and so $\a+ \beta$ or $\a-\beta$ is in $\Dl$.
Replacing $\beta$ by $-\beta$ if needed, we can assume that $\gamma := \a + \beta\in \Dl$.
Furthermore, $p(\gamma) = p(\a) + p(\beta) = p(\beta)\ne 0$ and so $\gamma\in \Dls$.
This proves 3.35.

Now $\Dl = \Dlns \cup (\Dl \setminus \Dlns)$ and, by 3.35, the sets
$\Dlns$ and $\Dl \setminus \Dlns$
are orthogonal.
But $\Dlns \supseteq \Dls \ne \emptyset$. Hence, by the irreducibility of $\Dl$
we have
$$\Dl =  \Dlns, \tag 3.36$$
which is statement (i).

To prove (ii), suppose  
that $$p(\Dls) = M_1 \cup M_2,\tag 3.37$$
where 
$M_1$ and $M_2$ are 
orthogonal (and hence disjoint).  
Since $p(\Dl)\setminus\{0\} = p(\Dls)$, it will suffice to prove that
$M_1$ or $M_2$ is empty.

Now by 3.37, we have the disjoint union
$$\Dls = K_1\cup  K_2,\tag 3.38$$
where $K_i = \{ \a\in \Dl \suchthat p(\a)\in M_i\}$, $i=1,2$.
We claim next that 
$$K_1 \text{ and }K_2 \text{ are orthogonal.}\tag 3.39$$
Indeed, suppose
for contradiction that there exist $\a\in K_1$, $\beta\in K_2$ with $(\a,\beta)\ne 0$.
Reasoning as above, we can assume that
$\a + \beta\in \Dl$.
But $(p(\a),p(\beta)) =0$ and so $p(\a+\beta) \ne 0$. Consequently,
$\a + \beta\in \Dls$ and therefore $\a + \beta$ is in $K_1$ or $K_2$.  We may
assume that $\a + \beta$ is in $K_1$.  Then, $p(\a + \beta)\in M_1$ and $p(\a)\in M_1$,
which implies that $(p(\beta),p(\beta)) = (p(\a+\beta) - p(\a),p(\beta) = 0$.
This contradiction proves 3.39.

Finally, 3.35, 3.36 and 3.38 tell us that $\Dl= \Psi_1 \cup \Psi_2$, where
$$\Psi_i = K_i \cup \{ \gamma - \beta \suchthat \beta, \gamma\in K_i,\ p(\gamma)= p(\beta),\ \gamma-\beta \in \Dl\},$$
$i=1,2$.
By 3.39, $\Psi_1$ and $\Psi_2$ are orthogonal.
Thus, $\Psi_1$ or $\Psi_2$ is empty, and so $K_1$ or $K_2$ is empty.  Hence,
$M_1$ or $M_2$ is empty as desired.
\qed\enddemo

We now return to our investigation of the axiom EA5a for $\tg$.
We first need to
recall some notation from Section 1. 
Let $V$ denote the real span of $R$ and 
assume that the form on $\g$ is scaled in such a way that
$(\a, \beta) \in \R $ for all $\a, \beta \in R$ and  the
form is positive semidefinite on $V$.
$V^0$ denotes the radical of $V$ and
the natural map from $V$ to $\bar V = V/V^0$ is denoted $x \mapsto \bar
x.$ Let $\bar R$ be the image of $R$ in $\bar V$.
Then $\bar R$ is a
finite irreducible root system in $\bar V$ where we use the
positive definite form on $\bar V$ induced from the form on $V$.
Notice that the map $x \to \bar x$  maps $R^\times$ onto the set
$\bar R^\times$ of nonzero roots of $\bar R$.

 From Lemma 3.23 we see that $\pi(R) \subseteq V$ and
hence $\pi(V) \subseteq V.$ In particular this implies that
$$V=(V \cap (\h^{\s})^*) \perp (V \cap ((\h^{\s})^{\perp})^*).
\tag 3.40$$
Furthermore, since the form is positive semidefinite on $V$ we get
that
$$V^0=(V^0 \cap (\h^{\s})^*) \perp (V^0 \cap ((\h^{\s})
^{\perp})^*). \tag 3.41 $$
It follows from this that we have, upon making the obvious
identifications, that
$$\bar V= 
\Big((V \cap (\h^{\s})^*)/(V^0 \cap (\h^{\s})^*)\Big)
\perp
\Big((V \cap ((\h^{\s})^{\perp})^*)/ (V^0 \cap ((\h^{\s})^{\perp})^*)\Big).
\tag 3.42$$
Note that 3.41 also implies that $\pi(V^0) \subseteq V^0,$ so that
$\pi$ induces a map $\bar \pi :\bar V \to \bar V.$ Moreover, it is
clear that $\bar \pi$ is just the projection onto the first factor in
3.42.

   We have seen that $\s(R) =R$ and so it follows that $\s(V)=V.$ Since
$\s$ preserves the form we get that
$\s(V^0)=V^0$ and so $\s$ induces an element $\bar \s$
of the orthogonal group of  our form on $\bar V.$
Furthermore, we have that
$$(V\cap (\h^{\s})^*)/(V^0 \cap(\h^{\s})^*) =\bar V ^{\bar
\s}, \tag 3.43 $$
and
$((V \cap ((\h^{\s})^{\perp})^*)/(V^0 \cap ((\h^{\s})^{\perp})
^*)$ equals the sum of the eigenspaces for $\bar \s$ corresponding
to  eigenvalues not equal to $1$. Also, it follows from Lemma 3.23
that
$$\bar \pi (\bar \a)=
\minv \sum_{i=0}^{m-1}\bar \s ^i (\bar \a)$$
for $ \bar \a \in \bar V.$

   We let $\tV$ be the real span of $\tR$ in $\th^*.$ Recall from 3.19
we have that
$$\tR= \bigcup_{i=0}^{m-1} (\pi(R_{\bari}) + i \Dl +m
\Z \Dl). \tag 3.44  $$
Since $0 \in R$ we have $0 \in R_{\bari}$ for some $i$ and so in
this case we obtain that $i \Dl +m \Z \Dl \subseteq \tR.$ Thus
$\Dl \in \tV$ so we obtain that $\pi(R_{\bari}) \subseteq
\tV$ for $0 \leq i \leq m-1.$ Thus $\tV$ is spanned by $\Dl$
together with $\pi(R)= \bigcup_{i=0}^{m-1}\pi(R_{\bar
\imath}).$ But by 3.40 the real span of
$\pi(R)$ is $V\cap(\h^{\s})^*.$ It follows that
$$\tV=(V \cap(\h^{\s})^*) \oplus \R \Dl= V^\s \oplus \R \Dl. \tag 3.45$$
Next we note that the restriction of the form from $\th^*$ to $V\cap
(\h^{\s})^*$ is the same as the restriction of the form from
$\h^*$ to $V\cap (\h^{\s})^*$.
It follows that the form on $\th^*$ restricted to $\tV$ is real valued
and positive semidefinite with $\Dl$ in its radical. Hence the
radical, $\tV^0$, of the form on $\tV$ is given by
$$\tV^0=(V^0\cap (\h^{\s})^*) \oplus \R \Dl =
(V^0)^\s \oplus \R \Dl.  \tag 3.46$$

Next we let $ \overline{\tV}= \tV / \tV^0$ and let $x \mapsto \bar x$ be
the associated natural homomorphism. Let $\overline {\tR}$ be the
image of $\tR$ in $\overline {\tR}.$ Then using 3.45 and 3.46 and
making the obvious identifications we have
$$\overline {\tV}= (V \cap(\h^{\s})^*)/(V^0 \cap (\h^{\s})^*)=
V^{\s}/(V^0)^{\s}. \tag 3.47 $$
In view of 3.42 and 3.43 we can identify $\overline \tV$ with the space
$\bar V ^{\bar \s}$ of $\bar \s$ fixed points  in $\bar V.$
Under this identification the element $\overline {\a_0+r \Dl}=
\a_0 +r \Dl +
\tV ^0 \in \overline {\tV}$ is identified with $\overline {\a_0}=
\a_0 +V^0 \in \bar V$ for $\a _0 \in V \cap (\h^{\s})^*$, $r
\in \R.$ Hence this identification preserves the forms on $\overline
{\tV}$ and $\bar V.$ But we have that $\overline { \tR}=\{ \overline
{\pi ( \a ) + i \Dl }\suchthat i \in \Z, \a \in R_{\bar
\imath} \},$ and
$\bar R = \{ \bar \a \suchthat i \in \Z, \a \in R_{\bar
\imath} \}.$ Hence under our identification
$$\overline {\tR}= \bar \pi (\bar R). \tag 3.48 $$

\proclaim {Lemma 3.49} We have
\roster
\item"{(i)}" If $\tR^{\times}$ is not empty and
$\a \in R^{\times}$,  then there exists
$\beta \in R^{\times}$ such that $(\pi (\beta), \pi( \beta ))\neq
0$ and $(\a, \beta) \neq 0.$
\item"{(ii)}" $\tR^{\times}$ cannot be written as the union
of two nonempty orthogonal subsets. Thus, EA5a holds for  $\tg.$
\endroster
\endproclaim

\demo{Proof}
To prove both of these statements we can assume that $\tR^{\times}$ is nonempty (since (ii)
is trivial if $\tR^{\times}$ is empty).  Then, $(\pi(\a),\pi(\a))\ne 0$ for some $\a\in R^\times$
by 3.30.  Thus, $\bar \pi(\bar \a) \ne 0$ for some $\bar \a \in \bar R^\times$.
So, since $\bar\pi  : \bar V \to \bar V$ is the orthogonal projection onto $\bar V^{\bar \s}$,
we see that $\bar V^{\bar \s}$ is nonzero.  Therefore, we can apply Lemma 3.34 to the
configuration $\Om = \bar R$,
$X = \bar V$, $Y = \bar V^{\bar \s}$ and $p = \bar \pi$.

To prove (i), let $\a\in R^\times$. Then, $\bar \a\in \bar R^\times$.  So, by Lemma
3.34~(i), there exists $\bar \beta\in \bar R^\times$ so that $\bar\pi(\bar\beta) \ne 0$
and $(\bar \a, \bar \beta) \ne 0$. Then, lifting $\bar\beta$ to an element $\beta$ of $R^\times$,
we have $(\pi(\beta),\pi(\beta)) =  (\bar\pi(\bar \beta),\bar\pi(\bar\beta)) \ne 0$ and $(\alpha,\beta) \ne 0$.

To prove (ii), suppose for contradiction that $\tR^\times = \tR_1\cup \tR_2$, where
$\tR_1$ and $\tR_2$ are orthogonal nonempty sets. 
Applying the map $x \mapsto \bar x$, we see that 
$\overline {\tR} \setminus \{ 0 \}$ is the union of two
nonempty orthogonal sets. 
By 3.48, the same is true
for $\bar\pi(\bar R) \setminus \{0\}$.  Since
$\bar\pi(\bar R) \setminus \{0\} = \bar\pi(\bar R^\times)\setminus \{0\}$,
this contradicts Lemma 3.34~(ii).
\qed\enddemo

   We next wish to compare the cores of our Lie algebras $\g$ and
$\tg$ (see Lemma 3.57 below).
Recall from [AABGP] that for $\a \in R^{\times}$ we have the
reflections $w_\a$ which satisfy
$$w_\a (\beta) =\beta - 2\frac{(\beta, \a)}{(\a, \a)}\a
\quad\text{for } \beta \in \h^*. \tag 3.50$$
Moreover, choosing $e_\a \in \g_\a, e_{-\a} \in
\g_{-\a}$ for $\a \in R^{\times}$ as in 1.18 of [AABGP] we
have the automorphisms (see 1.23 of [AABGP])
$$
\theta_\a(t) =\exp (\ad (te_\a)) \exp(\ad (-t^{-1}
e_{-\a}))\exp (\ad (te_\a)), $$
of $\g$ for any $t \in \C$ which satisfy
$$\theta_\a(t) \g_\beta =\g_{w_\a(\beta)} \text{ for any }
\beta \in R. \tag 3.51$$
We will let $\theta_\a$ denote the element $\theta_\a (1)$ in
what follows.

\proclaim{Lemma 3.52} Suppose that $\tR^{\times}$ is not empty. Then
the core, $\g_c$,  of  $\g$ is generated as an algebra by the root
spaces $\g_{\a}$ for $ \a \in R$, $(\pi (\a), \pi (
\a))\neq 0.$
\endproclaim

\demo {Proof}We let $\m$ be the subalgebra of $\g$ generated by all
the root spaces $\g_{\a}$ for $\a \in R, (\pi (\a), \pi
(\a )) \neq 0.$ Then $\m \subseteq \g_c.$ To prove the reverse
inclusion we  must show that $\m$ contains each root space $\g_\beta$
for $\beta \in R^{\times}.$ If $(\pi (\beta),\pi (\beta)) \neq 0$ this
is clear so we assume that
$(\pi (\beta ),\pi (\beta))=0.$ By Lemma 3.49~(i) there exists $\a \in
R^{\times}$ such that $(\pi (\a),\pi (\a)) \neq 0$ and
$(\a, \beta )\neq 0.$ Put $\epsilon = w_\a ( \beta) =\beta -2(
(\beta, \a)/(\a, \a))\a \in R^{\times}.$ Since $\pi
(\beta) \in V^0$ we have that $(\pi (\epsilon), \pi ( \epsilon)) =
(-2((\beta, \a)/(\a, \a)))^2(\pi (\a ), \pi (\a))
\neq 0.$ Then we have (using the notation developed before the
statement of our lemma) that $\theta_\a (\g_\epsilon)=
\g_{w_\a ( \epsilon)}=\g_\beta.$ But we know that
$\g_\epsilon$, $\g_\a$ and $\g_{ -\a}$ are contained in $\m.$
Hence, so is $\theta_\a (\g_\epsilon)$, and we obtain that $
\g_\beta \subseteq \m$ as desired.
\qed
\enddemo

\proclaim {Lemma 3.53} Suppose that $\tR^{\times}$ is not empty. Then
$\g_c$ is the sum of the spaces
$$\g_{\bari, \pi (\a)}\text { for } \bari \in \Z_m, \ \a \in
R_{\bari}, \ (\pi ( \a), \pi ( \a)) \neq 0, \tag 3.54 $$
together with their commutators.
\endproclaim

\demo{Proof} Let $\ss$ be the sum of the spaces in 3.54. We first show
that $\ss \subseteq \g_c.$ For this, let $\bari \in \Z_m$,
$\a \in \R_{\bari}$, $(\pi (\a),\pi (\a)) \neq 0.$
Then $\gia$ is contained in the $\pi ( \a)$-weight space for $\ad
(\h^\s).$ Thus $\gia \subseteq \sum g_\beta $ where the sum is
over those $\beta \in R$ satisfying $\pi (\beta) =\pi ( \a).$ But
any $\beta \in R$ satisfying $\pi ( \beta)= \pi ( \a)$ also
satisfies $(\pi (\beta), \pi ( \beta)) \neq 0$ and so is in $R^{
\times}.$ Hence $\ss \subseteq \g_c.$

   We next show that $\ss$ generates $\g_c$ as an algebra.  By Lemma
3.52 it is enough to show that $\ss$ contains all root spaces
$\g_\a$ with $\a \in R$, $(\pi (\a),\pi (\a)) \neq 0.$
But for such an $\a$, $\g_\a$ is contained in the $\pi ( \a)
$-weight space for $\ad (\h^\s)$ and so $\g_{\a}\subseteq
\sum_{i=0}^{m-1} \gia.$ Hence $\g_\a \subseteq \ss.$

   Since $\ss$ generates $\g_c$ as an algebra it only remains to show
that $\ss + [\ss, \ss]$ is closed under the Lie product and to show
this it is enough to show that
$$[\ss, [\ss, \ss]] \subseteq \ss + [\ss, \ss]. \tag 3.55 $$
For this let $i, j, k \in \Z$, $\a \in R_{\bari}$, $\beta \in
R_{\barj}$, $\epsilon \in R_{\bar k},$ with
$(\pi (\a), \pi (\a )) \neq 0, (\pi (\beta), \pi (\beta ))
\neq 0,(\pi (\epsilon), \pi (\epsilon )) \neq 0.$
We want to show that
$$[\gia,[\gjb, \gke] \subseteq \ss +[\ss,\ss ].$$
To prove this we can clearly assume that the left hand side is not
zero and so then $[\gjb,\gke]$ is also non-zero. Thus by 3.10, there
exists some $\eta \in R_{\barj +\bar k}$ so that $\pi ( \eta) =
\pi (\beta) +\pi (\epsilon)$ and
$[\gjb,\gke] \subseteq \g_{\barj + \bar k, \pi (\eta)}.$ So it
suffices to show that $$[\gia, \g_{\barj + \bar k, \pi (\eta)}]
\subseteq \ss + [\ss,\ss]. \tag 3.56$$

   If $(\pi ( \eta), \pi ( \eta)) \neq 0$ then the left hand side in
3.56 is in $[\ss,\ss].$ Thus, we can assume that
$(\pi ( \eta), \pi ( \eta)) = 0.$ But  there is some $\omega \in
R_{\bari +\barj + \bar k}$ such that
$ \pi ( \omega) = \pi ( \a) + \pi ( \eta)$ and
$$[\gia \g_{\barj + \bar k, \pi (\eta)}] \subseteq \g_{\bar
\imath +\barj +\bar k, \pi ( \omega)}.$$
Also, since $\pi ( \eta) \in V^0$ we have $ (\pi ( \omega), \pi (
\omega)) = (\pi ( \a), \pi ( \a )) \neq 0.$
It follows that we now have $[\gia,\g_{\barj + \bar k, \pi
(\eta)}] \subseteq \ss$ which is what we want.
\qed
\enddemo

Although $\tg$ may not be in general an EALA, we make the same definitions
of the core of $\tg$ and tameness of $\tg$ as in Section 1.
Thus, we define the \ital{core} of $\tg$ to be
the subalgebra $\tg _c$ 
of $\tg$ generated by $\tg_{\tilde\a}$, $\tilde\a\in \tR^\times$.  If
$\tR^\times$ is empty,
we take $\tg_c = \{0\}$.   We say that
$\tg$ is {\it tame\/} provided that the centralizer of $\tg_c$
in $\tg$ is contained $\tg_c$. In the next two lemmas, we calculate
the core of $\tg$ and give conditions which imply that $\tg$
is tame.

   The automorphism $\s$ of $\g$ stabilizes the core $\g_c$ of
$\g$ and hence it also stabilizes the center $Z(\g_c)$ of $\g_c$.
Thus we have
$$\g_c = \bigoplus_{i=0}^{m-1}(\g_c)_{\bar 
\imath}\andd Z(\g_c) = \bigoplus_{i=0}^{m-1} Z(\g_c)_{\bar 
\imath}$$
where  $(\g_c)_{\bari}=\g_{\bari} \cap \g_c$ and
$Z(\g_c)_{\bari}=\g_{\bari} \cap Z(\g_c)$ 
for $\bari \in \Z_m.$

\proclaim{Lemma 3.57} Suppose $\tR^{\times}$ is not
empty. Then
$$
\gather
\tg_c=\big(\bigoplus_{i \in \Z}(\g_c)_{\bari} \otimes   
t^i\big)
\oplus \C c, \tag 3.58\\
Z(\tg_c) = \big(\bigoplus_{i \in \Z}Z(\g_c)_{\bari} \otimes   
t^i\big)
\oplus \C c,\tag 3.59\\
\tg_c/Z(\tg_c) \cong \loopalg(\g_c/Z(\g_c),\s), \tag 3.60
\endgather
$$
where $\loopalg(\g_c/Z(\g_c),\s)$ denotes the loop algebra of
$\g_c/Z(\g_c)$ relative to the automorphism induced by $\s$
on $\g_c/Z(\g_c)$.
\endproclaim

\demo{Proof}Let $\ts$ denote the right hand side of 3.58. Since
$[(\g_c)_{\bari}, (\g_c)_{\barj}] \subseteq
(\g_c)_{\bari +\barj}$ for $ i,j \in \Z$ we obtain that
$\ts$ is a subalgebra of $\tg.$

   To prove that $\tg_c \subseteq \ts$ it suffices to show that
$\tg_{\ta} \subseteq \ts$ for $\ta \in \tR^{\times}.$
For this we let $\ta =\pi (\a) +i \delta$ where $i \in\Z$, $\a
\in R_{\bari}$, $(\pi (\a),\pi (\a))\neq 0.$ Then
$\tg_{\ta}= \gia \otimes   t^i$ which is contained in $(\g_c)_{\bar
\imath} \otimes   t^i$ since, by Lemma 3.53,
 $\gia$ is contained in  $\g_c.$ Thus, $\tg_c \subseteq \ts.$

   To prove that $\ts \subseteq \tg_c$ we first show that $c \in
\tg_c.$ Indeed, since $\tR^{\times}$ is not empty we can choose an
$\tilde \epsilon \in \tR^{\times}.$ Then $\tilde \epsilon =\pi
(\epsilon) + k \delta$ where $k \in \Z$, $\epsilon \in R_{\bar k}$,
$(\pi (\epsilon),\pi (\epsilon)) \neq 0.$ Then we also have that
$-\epsilon \in R_{-\bar k}$ and furthermore
$\g_{\bar k, \pi (\epsilon)}$ and $\g_{- \bar k ,\pi (-\epsilon)}$
are paired, in a non-degenerate manner, by the form $\fm$ on $\g.$
Thus, we can choose $x \in \g_{\bar k, \pi (\epsilon)}$ and $ y \in
\g_{-\bar k,\pi (-\epsilon)}$ so that
$(x,y) \neq 0.$ Also, $x \otimes t^k$ is in $\tg_{\tilde \epsilon}$
and $ y \otimes t^{-k}$ is in $\tg_{-\tilde \epsilon}$ and hence both
of these elements belong to $\tg_c.$ Similarly, since $x \in \g_{\bar
k +\bar m, \pi (\epsilon)}$ and $y \in \g_{-\bar k -\bar m, \pi(
-\epsilon)}$ we have $x \otimes t^{k +m} \in \tg_{\tilde \epsilon + m
\delta}$ and $y \otimes t^{-k -m} \in \tg_{-\tilde \epsilon -m\delta}$
and so the elements $x \otimes t^{k+m}$, $y \otimes t^{-k-m}$ are in
$\tg_c.$ Thus $[ x \otimes t^k, y \otimes t^{-k}] -[x \otimes t^{k +
m}, y \otimes ^{-k-m}]$ is in $\tg_c.$ But this element equals $-
m(x,y)c,$ and so we get that $c \in\tg_c$ as desired.

   For the proof of 3.58, it now only remains to show that $(\g_c)_{\bari} \otimes  
t^i \subseteq \tg_c$ for $ i \in \Z.$ But by Lemma 3.53,
$(\g_c)_{\bari}$ is the sum of spaces  of the form $\gia$,
$\a \in R_{\bari}$, $(\pi ( \a), \pi ( \a)) \neq 0,$
and $[\gjb, \g_{\bar l, \pi ( \mu)}],$ for $j,l \in \Z$, $\barj
+ \bar l =\bari$, $\beta \in R_{\barj}$, $\mu \in R_{\bar
l}$, $(\pi (\beta),\pi (\beta)) \neq 0$, $(\pi( \mu), \pi( \mu)) \neq
0.$ In the first case we have that $\gia \otimes   t^i= \tg _{\pi
(\a)+i \delta} \subseteq \tg_c.$ In the second case we can assume
that we have chosen $j$ so that $j+l=i.$ Then we have $[\gjb, \g_{\bar
l, \pi(\mu)}]\otimes   t^i \subseteq
[\gjb \otimes   t^j, \g_{\bar l, \pi (\mu)} \otimes   t^l] + \C c.$
But this is contained in
$[\tg_{\pi (\beta) + j \delta}, \tg_{\pi ( \mu) +l \delta}] +\C c
\subseteq [\tg_c,\tg_c] +\C c \subseteq \tg_c.$  So we have proved~3.58.

Since $\g_c$ is perfect, $Z(\g_c)$ is orthogonal to $\g_c$.  Thus,
we have the inclusion ``$\supseteq$'' in 3.59.  
The reverse inclusion in 3.59 follows from 3.58 and the fact that any element
of $Z(\tg_c)$ must commute with   
$(\g_c)_{\bari} \otimes t^j$ for all $j\in \Z$.
Finally, to prove 3.60, observe that by 3.58 and 3.59 we have
$$
\tg_c/Z(\tg_c)
\cong
\bigoplus_{i \in \Z}\big((\g_c)_{\bari}/Z(\g_c)_{\bari}\big) \otimes t^i
\cong 
\bigoplus_{i \in \Z}(\g_c/Z(\g_c))_{\bari} \otimes t^i = \loopalg(\g_c/Z(\g_c),\s)
$$
as vector spaces.  Moreover, it is clear  that this composite
map is an isomorphism of algebras.
\qed
\enddemo

We will actually not use 3.59 and 3.60 in what follows.  However,
the core modulo its center of an EALA plays a fundamental role in the
structure theory of EALA's (see [BGK], [BGKN] and [AG]).  Therefore
3.60 is of interest on its own.

\proclaim{Lemma 3.61} Suppose that $\g$ is tame and that $\tR^{\times}
$ is not empty. Then $\tg$ is also tame.
\endproclaim

\demo{Proof}Let $\tc$ be the centralizer of $\tg_c$ in $\tg.$ We must
show that $\tc \subseteq \tg_c.$

   Since $\tg_c$ is perfect, $\tc$ is orthogonal to $\tg_c.$ Thus,
$\tc$ is orthogonal to the central element $c$ by
 3.58. Hence $\tc \subseteq (\bigoplus_{ i \in\Z}\g_{\bari}
\otimes   t^i) \oplus \C c.$ 
Now let $x\in \tc$.  Then, $x =  \sum_{i\in \Z} x_i\ot t^i + rc$ for
some $x_i\in \g_{\bari}$, $r\in \C$. But $x$ commutes with $\tg_c$
and hence with $(\tg_c)_{\barj} \ot t^j$ for all $j\in \Z$ by 3.58.
Thus, for each $i$, $x_i$ commutes with $(\g_c)_\barj$ for all $j$
and so $[x_i,\g_c] = \{0\}$.  Hence, since $\g$ is tame, we have
$x_i\in (\g_c)_\bari$ and so $x\in \tg_c$ by~3.58.
\qed
\enddemo

   We now want to combine all of our previous work to conclude that
$\tg$ is a tame EALA provided that $\g$ is tame and $\tR^{\times}$ is not empty.
For this, we still need to see if EA5b holds for $\tg.$  But,
as we see in the next lemma, it is a general fact that if we have a
triple $\eala$ of two Lie algebras and a form which satisfies EA1,
EA2, EA3, EA4, EA5a and
if $\eala$ is tame, then $\eala$ automatically satisfies EA5b. No 
doubt, this result is of
independent interest. Of course, the core of such a triple, as well as
the concept of tameness, are defined in the usual manner.

\proclaim{Lemma 3.62} Let $\eala$ be a triple consisting of two Lie
algebras as well as a symmetric form and assume that EA1, EA2, EA3,
EA4 and EA5a hold. Further, let $\g_c$ be the subalgebra of $\g$
generated by the non-isotropic root spaces and assume that $C_{\g}
(\g_c)\subseteq \g_c.$ Then EA5b also holds and so $\eala$ is a tame
EALA.
\endproclaim

\demo{Proof} We use the notation from Section 1.
It is clear that $\g_c$ is perfect, that $C_{\g}(\g_c)$ is contained in
the orthogonal complement to $\g_c$, and that if $\delta$ is any root then
the restriction of the form to the space $\g_{\delta} +\g_{-\delta}$
is non-degenerate. In order to show EA5b holds we let $\delta$ be any
isotropic root of $\g.$

If $\g_{\delta} +\g_{-\delta} \subseteq C_{\g}(\g_c)$ then it follows,
since $C_{\g}(\g_c)$ is orthogonal to $\g_c$ and since $C_{\g}(\g_c)
\subseteq \g_c$, that $\g_{\delta} +\g_{-\delta}$ is totally
isotropic.  This is impossible and so $\g_{\delta} +\g_{-\delta}
\nsubseteq C_{\g}(\g_c)$. Thus, since the non-isotropic root spaces
generate $\g_c$, it follows that there exists a non-isotropic root $\a
\in R$ such that either $[\g_{\a},\g_{\delta}] \neq \{ 0 \}$ or
$[\g_{\a}, \g_{-\delta}] \neq \{ 0 \}.$ If
$[\g_{\a},\g_{\delta}] \neq \{ 0 \}$ then we get that $\a +
\delta \in R$. If we have $[\g_{\a}, \g_{-\delta}] \neq \{ 0 \}$
then we get that $-\delta + \a \in R$ and so
its negative $\delta -\a$ is in  $R.$ In either case we get what we want.
\qed
\enddemo

Combining the results of this section, we have the following theorem.

\proclaim {Theorem 3.63} Let $\eala$ be a tame extended affine Lie
algebra, let $m$ be a positive integer and let $\s$ be an
automorphism of $\g$ satisfying
\smallskip
{\bf A1.} $\s^m=1$

{\bf A2.} $\s(\h)=\h$

{\bf A3.} $(\s(x),\s(y))=(x,y)$ for all $x,y \in \g$

{\bf A4.} The centralizer of $\h^{\s}$ in $\g^{\s}$ equals $\h^{\s}$.

\smallskip\noindent
Let $\tg= \affgs$, $\th=\h^{\s} \oplus \C c \oplus \C d$ and let
$\fm$ be the form defined by 2.4 (restricted to $\affgs$). Let $\tR^\times$
be the set of nonisotropic roots $\tg$ relative to $\th$ .
Then either $(\tg,\th, \fm)$ is a tame EALA  
or $\tR^{\times}$ is empty.
Furthermore, $\tR^{\times}$ is nonempty 
if and only if
$$(\pi ( \a), \pi( \a))\ne 0 \text{ for some } \a \in R, 
\tag
3.64$$
where $\pi ( \a)= \minv\sum_{i=0}^{m-1}\s^i(
\a).$
\endproclaim

\demo{Proof}
The first statement follows from
3.28, 3.32, 3.33, 3.49~(i), 3.61 and 3.62.  The second  statement follows
from  3.30.\qed\enddemo

\medskip

If $\eala$ is a tame EALA, $\s$ satisfies A1, A2 and A3, and $m$ 
is a prime, the following corollary tells us that we can
test whether or not $(\tg,\th, \fm)$ is a tame EALA using only
information about the root system $R$ of $\g$
and the action of $\s$ on $R$.  This corollary follows
immediately from the theorem and Proposition 3.25.

\proclaim {Corollary 3.65} Let $\eala$ be a tame EALA, let $m$ be a 
positive integer and let $\s$ be an
automorphism of $\g$ satisfying A1, A2 and A3.
Let $(\tg,\th, \fm)$ be as in Theorem 3.63.  
If $\pi(\a) \ne 0$ for all nonzero $\a\in R$ and $(\pi(\a),\pi(\a)) \ne 0$
for some $\a\in R$, then 
$(\tg,\th, \fm)$ is
a tame EALA.  Moreover, the converse is true if $m$ is a prime.
\endproclaim

\remark{\rm \bf Remark 3.66}  Let $\eala$ be a tame EALA, let 
$\s$ be an automorphism of $\g$ which 
satisfies A1, A2, A3 and A4, and suppose that $\Rt^\times$ is 
nonempty.  Then $(\tg,\th, \fm)$ is a tame EALA.  
It is interesting to compare the type and nullity of $\g$ with the 
type and nullity of $\tg$.  Recall that the type of $\g$ is the type 
of the finite root system $\bar R$, whereas the nullity $\nu$ of $\g$ 
is $\dim_\R V^0$, where $V^0$ is the radical of the real span $V$ of 
the roots of $\g$. 

\smallskip
(i) One can compute the type of $\tg$  by calculating the finite root 
system $\overline {\tR}=\bar \pi ( \bar R)$ using the formula $\bar 
\pi ( \bar \a )=\minv \sum_{i=0}^{m-1}\bar \s 
^i(\bar \a)$ (see 3.48).

\smallskip
(ii) The nullity $\tilde \nu$ of $\gt$ is given by
$$\tilde \nu = \dim_\R((V^0)^\s) + 1 \le \nu +1,$$
where $(V^0)^\s$ is the space of fixed points of $\s$ acting 
on $V^0$
(see 3.41 and 3.45).
\endremark

\proclaim {Corollary 3.67} Let $\eala$ be a tame EALA. Let $\affg$ be the affinization of $\g$
and let $\fm$ be the form defined by 2.4.
Then, the triple 
$$(\affg,\h\oplus \C c \oplus\C d, \fm)$$ is a tame EALA.   Moreover, the type of $\affg$ is the same as 
the type of $\g$, and the nullity of $\affg$ is one more than the 
nullity of  $\g.$
\endproclaim

\demo{Proof}  
Let $\s
=\ident_{\g}$ and $m=1$. 
Clearly A1, A2, A3 and A4 hold.  Also $\pi(\a) = 
\a$ for
$\a\in R$ and so 3.64 is clear.  Thus, by Theorem 3.63, 
$\tg = \affg$ is a tame EALA.
Furthermore, by Remark 3.66~(i), $\overline {\tR}= \bar R$ and so 
$\tg$ and $\g$ have the same type.  Also, by Remark 3.66~(ii), we 
have $\tilde \nu = \nu +1$.
\qed
\enddemo

We conclude this section with a brief discussion of degeneracy of EALA's.
Recall from [G] that an EALA $\eala$ is said to be \ital{non-degenerate} 
if the dimension of the real vector space  $V^0$ spanned by the isotropic roots $R^0$ of $\g$ equals the dimension of the complex vector space spanned by $R^0$. 
The assumption that an EALA is nondegenerate simplifies the description of its structure.  Thus, it is of interest to describe the extent to which affinization preserves degeneracy.

If $W$ is a subset of a complex vector space, we  let $W_{\C}$ denote the complex space spanned by $W$. Using this notation, we have that
$$ \g \text{ is non-degenerate if and only if } \dim_{\C}((V^0)_{\C})=\dim_{\R}(V^0). $$
Otherwise said, $\g$ is non-degenerate if and only if $\dim_{\C} ((V^0)_{\C})=\dim_{\C}(V^0 \otimes_{\R}\C).$

We now assume that we have an EALA $\eala$ and an automorphism  $\s$ of $\g$
which satisfies A1, A2, A3, and A4.  Let $\tg$ and $\th$ be as in Theorem 3.63. 
Although $\tg$ is not in general an EALA, we nevertheless say that $\tg$ is
\ital{non-degenerate} if $\dim _{\C}((\tV^0)_{\C})=\dim_{\R}(\tV^0).$

\proclaim {Proposition 3.68}  Let $\eala$ be a tame non-degenerate extended affine Lie algebra and let $\s$ be an automorphism of $\g$ which satisfies A1, A2, A3 and A4. Let $(\tg,\th,\fm)$ be as in Theorem 3.63. Then $\tg$ is non-degenerate.  
In particular, if $\tg$ is an EALA then $\tg$ is a non-degenerate EALA.
\endproclaim

\demo{Proof}
We are assuming that 
$$\dim_{\C}((V^0)_{\C})=\dim_{\R}(V^0),\tag 3.69$$
and we want to prove that
$\dim _{\C}((\tV^0)_{\C})=\dim_{\R}(\tV^0).$ 
Now by 3.46 we have that $\tV^0=(V^0)^\s \oplus \R \delta$
and so $\dim_{\R}(\tV^0)=\dim_{\R}((V^0)^\s)+1$. 
But we also have 
$((V^0)^\s)_\C \cap \C\delta = \{0\}$ by the very definition of $\delta$  (see 3.16 and 3.17).  Hence,
$(\tV^0)_{\C} 
= ((V^0)^\s)_\C \oplus \C \delta$ and so
$\dim_\C((\tV^0)_\C)= \dim_{\C}(((V^0)^\s)_\C)+1$. 
Thus it is enough to show that 
$$\dim_{\C} (((V^0)^\s)_{\C})=\dim_{\R}((V^0)^\s). \tag 3.70$$

Now 3.69 is equivalent to the statement that the natural $\C$-linear map
$V^0\otimes_\R \C \to (V^0)_\C$ is an isomorphism.   Equivalently, 3.69 
says that every subset of $V^0$ that is independent over $\R$ is independent over $\C$.
This property is inherited by
real subspaces of $V^0$ and so 3.70 holds.
\qed
\enddemo

\head \S 4 Examples \endhead

In this section, we describe three examples that illustrate the use of
our main theorem 
to construct an EALA 
$\affgs$ starting from an EALA $\g$ and an automorphism $\s$ of finite order.

In the first example, $\g$ is a toroidal Lie algebra
coordinatized by the ring of commutative Laurent polynomials in $\nu$ variables.

\example{Example 4.1} Let 
$$\A = \Bbb{C}[t_1^{\pm 1},\cdots, t_\nu^{\pm 1}]$$
be the ring of
commutative Laurent polynomials in $\nu$ variables.  In $\A$ we use
the standard notation $t^{\p} = t_1^{p_1}\dots t_\nu^{p_\nu}$ for
$\p = (p_1,\dots,p_\nu) \in \zn$. Thus $\A$ is a
$\zn$-graded commutative associative algebra with $\A^{\p} = \Bbb{C}t^{\p}$
for $\p\in\zn$.  We will make use of the 
$\C$-linear map
$\epsilon : \Cal{A} \to \Bbb{C}$ given by linear extension of
$$\epsilon (t^{\p}) =
\cases  1 &\text{ if } \p = 0\\
0 & \text{  if  } \p \neq 0\endcases.$$

Let $\dot \g$ be a finite dimensional simple
complex Lie algebra. Let $\tau$ be an automorphism of $\dg$ of period
$m.$  We fix a 
Cartan subalgebra $\k$ of the 
nonzero reductive Lie algebra $\dg^{\tau}$ (see [BM] and [P]). Then the 
centralizer $\dh$ of $\k$ in $\dg$ is a Cartan subalgebra of $\dg$ 
which is stable under $\tau$ (ibid). Clearly $\dh^{\t} = \k.$

Based on this choice of $\dg$ and $\dh$, we construct a tame EALA $\eala$ along 
the lines of Example 1.29 of Chapter 3 of [AABGP]. We recall how this goes.

First let $\K = \dg \otimes \A$. Then, $\K$ is a $\Z^\nu$-graded
Lie algebra with $\K^\p = \dg \otimes \A^\p$.  
Let $\fm_{\K}$ be the unique bilinear form 
on $\K$
satisfying $(\dot x \otimes a,\dot y \otimes b)_{\K} = (\dot x,\dot y)\epsilon (ab)$ for 
$\dot x,\dot y\in\dot \g$, $a,b\in\A$, where $(\dot x,\dot y)$ is  the Killing form on $\dot \g.$ 
Then $\fm_{\K}$ is a nondegenerate invariant symmetric bilinear form on $\K$.

As vector spaces we let
$$\g = \K \oplus\Cal C \oplus \Cal D,$$
where $\Cal{C}= \Bbb{C}c_1
\oplus \cdots \oplus \Bbb{C}c_\nu$
and $\Cal{D} = \Bbb{C}d_1 \oplus
\cdots \oplus \Bbb{C} d_\nu$
are $\nu$-dimensional. 
The bracket on $\g$ is defined so that
$$\gathered
[\g, \Cal{C}] = [\Cal{D}, \Cal{D}] = \{0\},\\
[d_i, x ] = p_i x \ \text{  for all } 1 \leq i \leq \nu \text{ and }x\in \K^\p,\\
[x, y] = [x, y]_{\K} + \sum_{i=1}^\nu ([d_i, x], y)_\K
c_i \text{ for all } x, y \in \K,
\endgathered \tag 4.2$$
where $[ \cdot ,\cdot]_{\K}$ denotes the bracket on $\K$.
The bilinear form $\fm$ on $\g$ is defined so that
$$\gathered
\text{$\fm$ extends $\fm_{\K}$},\\
(\Cal{C}, \Cal{C}) = (\Cal{D}, \Cal{D})= 
(\Cal{C}, \K) = (\Cal{D},\K) = \{0\} \quad\text{and}\\ 
(c_i, d_j) = \delta_{ij}, \ i, j = 1,\dots, \nu.
\endgathered\tag 4.3$$
Finally, we set
$$\h = \dh  \oplus\Cal C \oplus \Cal D$$
where we are identifying $\dh = \dh \otimes 1$.
Then $\eala$ is a tame EALA.

We now describe the automorphism $\s$ of $\g$ that we will use.
First we fix $\boldmu = (\mu_1,\dots,\mu_\nu) \in \zn.$  
Then $\boldmu$ 
uniquely defines an 
automorphism  $\mu \in Aut(\A)$ satisfying 
$\mu(t^{\p}) = \zeta^{\boldmu \cdot \p}t^{\p}$ where $\boldmu 
\cdot  
\p  := \mu_{1}p_{1} + \ldots + \mu_{\nu}p_{\nu}.$  
Both $\mu$ and $\tau$ admit natural extensions (that again we will 
denote by $\mu$ and $\tau$) to $\Aut(\g)$ as follows:
$$\mu \text{  acts like } 1 
\otimes \mu \text{ on } \dg \otimes \A, \text{ and } \mu 
\text { fixes } \Cal C \text{ and  } \Cal D 
\text{ pointwise} $$
and
$$\tau \text{  acts like } \tau 
\otimes 1 \text{ on } \dg \otimes \A, \text{ and } \tau 
\text{ fixes } \Cal C \text{ and  } \Cal D 
\text{ pointwise.} $$
We set
$$\s := \t \mu.$$

We claim that $\s$ satisfies conditions A1 thru A4. 
Indeed A1 
is clear since $\t$ and $\mu$ commute. A2 holds because both $\t$ 
and $\mu$ stabilize $\dh$, $\Cal C$ and $\Cal D.$ A3
follows from the fact that the Killing form on $\dg$
is 
$\tau$-invariant. 
To check A4, first note that
$$\g^{\s} =(\bigoplus
\Sb
\p \in \zn \\
\endSb
\dg_{-\overline{\boldmu \cdot \p}} \otimes \Cal{A}^{\p})\oplus\Cal C \oplus \Cal D
\andd
\h^{\s} = (\k \otimes 1) \oplus\Cal C \oplus \Cal D $$
If $\p \in \zn \setminus {0}$ and $x \in \dg,$ then $[\Cal D, x 
\otimes \Cal{A}^{\p}] = \{0\}$ forces $x = 0.$ As a consequence
$C_{\dg^{\s}}(\h^{\s})$ is contained in  
$(\g_{\overline 0}\ot \A^{0})\oplus \Cal C \oplus \Cal D = 
(\dg^{\t}\ot 1) \oplus \Cal C \oplus 
\Cal D.$
But $C_{\dg^{\t}}(\k) = \k$ (since $\k$ is a Cartan subalgebra of 
$\dg^{\t}$), and therefore $C_{\g^{\s}}(\h^{\s}) = 
\h^{\s}$ as prescribed by A4.

We now consider the triple $(\tg,\th,\fm)$, 
where
$\tg=\affgs= (\bigoplus_{i\in \Z} \g_{\bari}\otimes t^i)
\oplus \C c \oplus \C d,$
$\th = \h^\sigma \oplus \C c \oplus \C d$
and $\fm$ is the restriction of the form 2.4.
We want to next show that $\Rt^\times$ is nonempty.
By Theorem 3.63,
it will suffice to show that $(\pi(\a), 
\pi(\a)) \neq 0$ for some $\a \in R.$   To see this, we identify $\dh^*$ as a subspace
of $\h^*$ in the obvious way.  Then the finite root system $\dot R$ of
$\dg$ with respect to $\dh$ becomes identified with a subset 
of $R$.  Now
since $\{0 \} \neq \k \subset \dg$ 
and $\dot R$ spans $\dh^{*},$ there exists an $ \a \in \dot R$ such 
that $\a \mid _{\k} \neq 0$. 
But then $\pi(\a) \ne 0$ by definition of $\pi$. So
$(\pi(\a), \pi(\a)) \ne 0$, since the Killing form restricted to the real span of $\dot R$
is definite.  Thus, by Theorem 3.63, $\Rt$ is nonempty and hence $(\tg,\th,\fm)$ is a
tame EALA.

It is interesting to note that the isomorphism class of the EALA
$\tg = \Aff(\g,\s)$ just constructed depends on $\tau$ but not on $\mu$.  This fact is a consequence of an ``erasing'' result
that we will prove in the next paper in this series.
\endexample

\bigskip
In our next example, $\g$ is an EALA of type $A_\ell$ with a noncommutative
coordinate algebra.

\example{Example 4.4}  
Let $\q = (q_{ij})$ be a $\nu\times\nu$ complex matrix so that
$q_{ii} = 1$ and $q_{ij} = q_{ji}^{-1}$.  
Let 
$$\A = \Bbb{C}_\q[t_1^{\pm 1},\cdots, t_\nu^{\pm 1}]$$
be the \ital{quantum torus} determined by $\q$.  Thus, by definition,
$\A$ is the associative algebra
generated by $t_1^{\pm 1},\dots,t_\nu^{\pm 1}$
subject to the relations $t_i {t_i}^{-1} = {t_i}^{-1}t_i = 1$ and
$t_it_j = q_{ij}t_jt_i$.
In $\A$ we write $t^{\p} = t_1^{p_1}\dots t_\nu^{p_\nu}$ for
$\p = (p_1,\dots,p_\nu) \in \zn$, in which case
$\A$ is a $\zn$-graded associative algebra with $\A^{\p} = \Bbb{C}t^{\p}$
for $\p\in\zn$.  We define $\epsilon : \Cal{A} \to \Bbb{C}$ exactly as in 
the commutative case (Example 4.1).

The triple $\eala$ that we use was introduced in [BGK].
We briefly recall the description of this EALA. 

First of all let $\K = \sll_{\ell+1}(\A)$.   That is, $\K$ is the Lie algebra
of $(\ell+1)\times(\ell+1)$ matrices over $\A$ generated by
the elementary matrices $ae_{i,j}$, $1\le i\ne j\le \ell+1$, $a\in\A$.
$\K$ has a unique $\Z^\nu$ grading so that $ae_{i,j}\in \K^\p$ for 
$1\le i\ne j\le \ell+1$, $a\in\A^\p$.  We define a nondegenerate invariant
symmetric  bilinear
form $\fm_{\K}$ on $\K$ by setting $(x,y)_\K = \epsilon(\tr(xy))$ for
$x,y\in \K$. 
As vector spaces we let
$$\g = \K \oplus\Cal C \oplus \Cal D,$$
where $\Cal{C}= \Bbb{C}c_1
\oplus \cdots \oplus \Bbb{C}c_\nu$
and $\Cal{D} = \Bbb{C}d_1 \oplus
\cdots \oplus \Bbb{C} d_\nu$
are $\nu$-dimensional. 
As in Example 4.1, we define a bracket $[\, , \,]$ and a form
$\fm$ on $\g$ by (4.2) and (4.3).
Finally, we set
$$\h = \dh  \oplus\Cal C \oplus \Cal D,$$
where $\dh = \sum_{i=1}^{\ell} \C(e_{i,i}-e_{i+1,i+1})$.
Then $\eala$ is a tame EALA.  The set $R$ of roots of $\g$ is
given by 
$$R = \{\sum_{k=1}^\nu n_k \delta_k \suchthat
\ n_i\in \Z\} \cup
\{\e_i - \e_j +\sum_{k=1}^\nu n_k \delta_k \suchthat
1\le i \ne j \le \ell+1,\ n_i\in \Z\},\tag 4.5
$$
where $\e_i(h)$ is the $i^{\text{th}}$ entry of $h$ for $h\in \dh$,
$\e_i$ is zero on $\Cal C$ and $\Cal D$, 
$\delta_i(d_j) = \delta_{i,j}$ and
$\delta_i$ is zero on $\dh$ and $\Cal C$.

We now assume (in addition to our earlier assumptions) that $q_{ij} = \pm 1$ for all $i,j$.
Then there exists an involution
(antiautomorphism of period 2) \ $\inv$ \ of $\A$ so that $\overline{t_i}= t_i$
for all $i$ (see for example [AG, \S 2]).  Since $\overline{t_1^{p_1}\dots t_\nu^{p_\nu}} = t_\nu^{p_\nu}\dots t_1^{p_1}$, 
the map \ $\inv$ \ is called the \ital{reversal involution} on $\A$.
Using \ $\inv$ \ we can define an involution $*$ on the associative algebra
$M_{\ell+1}(\A)$ by  $(a_{ij})^* = (\overline{a_{\ell+2-j,\ell+2-i}})$.  
Then, we define  a linear map $\s : \g \to \g$  by setting
$$\gather
\s(x) = -x^* \quad\text{for } x\in \K,\\
\s|_{\Cal C} = \id_{\Cal C} \andd \s|_{\Cal D} = \id_{\Cal D}.
\endgather$$
Observe that  if $x,y\in \K$, we  have 
$\epsilon\big(\tr(yx)\big)=\epsilon\big(\tr(xy)\big)$ (see p.~366 of [BGK]),
and so $(\s x, \s y)_\K = \epsilon\big(\tr((\s x)(\s y))\big) =
\epsilon\big(\tr(x^*y^*)\big) = \epsilon\big(\tr((yx)^*)\big)=\epsilon\big(\tr(yx)\big)=\epsilon\big(\tr(xy)\big) = (x,y)_\K$.
Using this fact is easy
to check that $\s$ is an automorphism of $\g$ that preserves the form $\fm$.
In fact it is clear that $\s$ satisfies axioms A1, A2 and A3, where $m=2$.

We next check A4 using Proposition 3.25.  Indeed, one has
$$\gathered
\s\big(\sum_{k=1}^\nu n_k \delta_k\big) = \sum_{k=1}^\nu n_k \delta_k 
\quad\text{and}\\
\s\big(\e_i - \e_j +\sum_{k=1}^\nu n_k \delta_k\big) 
=  \e_{\ell+2-j} - \e_{\ell+2-i}  +\sum_{k=1}^\nu n_k \delta_k.
\endgathered \tag 4.6
$$
Thus, if $\a$ is a nonzero root in $R$, we have
$\s(\a) \ne -\a$  and so $\pi(\a) \ne 0$.  So A4 holds.  

Once again we consider the triple $(\tg,\th,\fm)$ as in Theorem 3.63.  Observe
that $\pi(\e_1-\e_2) = \frac 12 (\e_1-\e_2+\e_\ell-\e_{\ell+1})$ is nonisotropic.  Hence, by Theorem 3.63,
$\Rt$ is not empty and so $(\tg,\th,\fm)$ is a tame EALA.

Using 4.5, 4.6 and Remark 3.66, we can calculate the type and nullity of $\tg = \affgs$.  
Indeed, we have 
$$\bar R =\{\bar 0\} \cup
\{\bar\e_i - \bar\e_j \suchthat
1\le i \ne j \le \ell+1\}.$$
Thus, the finite root system
$\overline{\Rt}$ associated with $\tg$ is
$$\overline{\Rt} = \bar\pi(\bar R) = \{\bar 0\} \cup
\{\frac 12(\bar\e_i - \bar\e_j  + \bar\e_{\ell+2-j} - \bar\e_{\ell+2-i})\suchthat
1\le i \ne j \le \ell+1\}.$$   One easily sees that the type of this
finite root system is $C_p$ if $\ell = 2p-1$ and $BC_p$ if $\ell = 2p$.
Hence $\tg$ has type $C_p$ if $\ell = 2p-1$ and type $BC_p$ if $\ell = 2p$.
Also, since all isotropic roots are fixed by $\s$, the nullity of $\tg$ is $\nu + 1$.

This example can be regarded as a noncommutative and higher nullity version of Kac's construction
of the affine Lie algebra of type $A_{\ell}^{(2)}$ from the finite dimensional
simple Lie algebra of type $A_{\ell}$.  (In fact this affine
Lie algebra is alternately denoted by $C_\p^{(2)}$
if $\ell = 2p-1$ and $BC_p^{(2)}$ if $\ell = 2p$. See [MP].)

\endexample

\bigskip
In the next example, we consider the case when the Lie algebra $\g$ is an affine Kac-Moody Lie algebra and the automorphism $\s$ is a diagram automorphism.

\example{Example 4.7} Let $\g = \g(A) = \g'\oplus \C d$ be the affine Kac--Moody Lie  
algebra
constructed from an affine $(\ell+1)\times(\ell+1)$ generalized Cartan matrix $A$.  
We use the notation of [K2, Chapter 6].  Let $\fm$
be the normalized standard invariant form on $\g$,  and let
$\h$ be the Cartan subalgebra of $\g$ used in the definition of $\g$. Then
$\eala$ is a tame EALA of nullity 1.

Let $\Pi =
\{\a_0,\dots,\a_\ell\}$ be the root basis of $\g$
and let
$a_0,\dots,a_\ell$ denote the
unique relatively prime positive integers so that
$A [a_0, \hdots, a_\ell]^t = 0$.
Then, the lattice $R^0$ of isotropic roots is given by $R^0 = \Z \delta$,
where
$\delta = \sum_{i=0}^\ell a_i  \a_i$.
Consequently the real span $V^0$ of $R^0$ is given by $V^0 = \R \delta$.

Next let $\s$ be an
automorphism of period $m$ of the GCM $A$.  Thus, by definition, $\s$ is a 
permutation of period $m$ of the set $\{0,...,\ell\}$
and $\s$ satisfies $a_{\s i,\s j} = a_{i,j}$
for all $i,j$.
We recall how to  ``extend'' $\s$ to an automorphism of $\g$. 
Indeed it is shown in [FSS, \S 3.2] that there is a unique automorphism of $\g$,
which we also denote by $\s$, so that $\s$ preserves the form $\fm$ on $\g$ and
$$\s e_i = e_{\s i} \andd \s f_i = f_{\s i}$$
for all $i$.  Moreover,  $\s$ stabilizes $\h$
and has period $m$ (ibid). We call the automorphism $\s$ of $\g$
the {\it diagram automorphism} associated with
the automorphism $\s$ of the GCM~$A$. 

Now $\s$ acts on the set $R$ of roots of $\g$.
Moreover since $\s e_i = e_{\s i}$, we have
$$\s\a_i = \a_{\s i}$$
for all $i$. Thus $\sigma$ permutes
the elements of $\Pi$. 

It is clear that the diagram automorphism $\s$ satisfies A1, A2 and A3.
We now check A4.  By Proposition 3.25, it is enough to show that $\pi(\a)\ne 0$
for all nonzero $\a$ in $R$.
For this we can assume that $\a$ is positive.  
But $\s$ permutes $\Pi$ and stabilizes $R$.  Thus
$\s^i(\a)$ is a positive root for $i\ge 0$.
Hence $\pi(\a) =  \frac 1m \sum_{i=0}^{m-1} \s^i \a$ is nonzero.
So A4 holds.

We now consider the triple $(\tg,\th,\fm)$, where 
$\tg=\affgs= (\bigoplus_{i\in \Z} \g_{\bari}\otimes t^i)
\oplus \C \tilde c \oplus \C \tilde d,$
$\th = \h^\sigma \oplus \C \tilde c \oplus \C \tilde d$
and $\fm$ is the restriction of the form 2.4.
(Here we use the notation $\tilde c$, $\tilde d$, and $\tilde \delta$ 
for the Lie algebra $\gt$, since we
have already used $c$, $d$ and $\delta$ for~$\g$.)  
We next want to decide when the set $\Rt^\times$ of nonisotropic roots of $\tg$
is non empty.  For this purpose,
we consider cases.

Suppose first that $\sigma$ (or more precisely the group
generated by $\s$) acts transitively on $\Pi$.   Now from the uniqueness
of the integers $a_i$ it follows that $a_{\s i} = a_i$ for all $i$.
Thus, by transitivity, all of the integers $a_i$ are equal and so
they are all $1$.  Thus, for each $j$, we have
$$\pi(\a_j) = \frac 1m \sum_{i=0}^{m-1} \s^i \a_j
= \frac 1m \frac m{\ell+1}\sum_{i=0}^{\ell} \a_i= \frac 1{\ell+1}\sum_{i=0}^{\ell} \a_i= \frac 1{\ell+1} \delta.$$
Consequently, $\pi(\a_j)$ is isotropic for all $j$, and so $\pi(\a)$ is isotropic
for all $\a\in R$.  Hence, by Theorem 3.63, $\Rt^\times$ is empty.

Suppose next that $\sigma$ does not act transitively on $\Pi$. 
Fix $j\in \{0,\dots,\ell\}$. Then
since $\s$ does not act transitively,  
$\pi(\a_j)$ lies in the real span of a proper subset of $\Pi$. 
But, as we saw when verifying A4, $\pi(\a_j)\ne 0$.
Thus $\pi(\a_j)$ is not isotropic.  So by Theorem 3.63, 
$\Rt^\times$ is non empty and $(\tg,\th,\fm)$ is a tame EALA. 
Finally, to calculate the nullity of $\tg$, notice that $\s \delta$ is a positive
root that generates $R^0$.  Hence, $\s \delta = \delta$ and so
$(V^0)^\sigma = V^0$ has real dimension 1.  Thus, by Remark 3.66~(ii),
$\gt$ has nullity 2.
\endexample

We summarize the conclusions from this example in the following theorem.

\proclaim{Theorem 4.8} Suppose that $\s$ is a diagram 
automorphism of
an affine Kac-Moody Lie algebra $\g = \g(A)$.   
Then $\s$ satisfies axioms A1, A2, A3 and A4.  Furthermore,
let $\tg=\affgs$ and $\th = \h^\sigma \oplus \C \tilde c \oplus \C \tilde d$,
and let $\fm$ be the restriction of the form 2.4.  Then
\roster
\item"{(i)}"  If $\sigma$ acts transitively on $\Pi$, the 
set $\Rt^\times$ of nonisotropic roots of $\tg$ is empty.
\item"{(ii)}"  If $\s$ does not act transitively on $\Pi$, $(\tg,\th,\fm)$ is  a tame EALA of nullity~2.
\endroster
\endproclaim

We note that case (i) in the theorem 
can only occur for one affine GCM, namely $A =
A_\ell^{(1)}$, $\ell \ge 1$.  (This follows from the
classification of affine GCM's.) Furthermore, if
$A =A_\ell^{(1)}$ and we label the roots
of $\Pi$ as in [K2, \S 6.1], then $\sigma$ acts transitively on $\Pi$
if and only if 
$\s = \tau^t$,
where $\tau = (0,1,\dots,\ell)$ is the diagram rotation and $t$ is
relatively prime to $\ell+1$.

\enddocument